\documentclass{article}
\usepackage{authblk}
\usepackage{xcolor}
\usepackage{amsthm}
\usepackage{amsmath}
\usepackage[utf8]{inputenc}
\usepackage{amssymb}
\usepackage{tikz}
\usepackage{geometry}
\usepackage{todonotes}
\usepackage{graphics}
\usepackage{thm-restate}
\usepackage{authblk}
\usepackage{multicol}
\usepackage[T1]{fontenc}
\usepackage[english]{babel}
\usepackage{graphicx}
\usepackage{fullpage}
\usepackage{caption}
\usepackage{subcaption}
\usepackage{cleveref}
\usepackage{lipsum}

\usetikzlibrary{calc}
\newcommand{\decisionpb}[4]{
        \begin{minipage}{#4\textwidth}
                #1\\
                \emph{Instance:} #2\\ 
                \emph{Output:} #3
        \end{minipage}
}

\newtheorem{theorem}{Theorem}
\newtheorem{remark}[theorem]{Remark}
\newtheorem{lemma}[theorem]{Lemma}
\newtheorem{corollary}[theorem]{Corollary}

\newtheorem{definition}[theorem]{Definition}

\newtheorem{conjecture}[theorem]{Conjecture}

\newenvironment{claim}[1]{\par\noindent\underline{Claim:}\space#1}{}
\newenvironment{claimproof}[1]{\par\noindent\underline{Proof:}\space#1}{\hfill $\blacksquare$}

\newcommand{\inc}{\textsc {Incidence}}
\newcommand{\mminc}{\textsc {Maker-Maker Incidence}}
\newcommand{\mbinc}{\textsc {Maker-Breaker Incidence}}

\newcommand{\MBPG}{\textsc{Maker-Breaker Positional Game}}
\newcommand{\MBSPG}{\textsc{Maker-Breaker Scoring Positional Game}}
\newcommand{\MMSPG}{\textsc{Maker-Maker Scoring Positional Game}}

\tikzstyle{v}=[circle, inner sep=0, minimum size =6 pt, line width = 1pt, draw=black, fill=white, text= black]
\tikzstyle{L}=[circle, inner sep=0, minimum size =6 pt, line width = 1pt, draw=black, fill=blue, text= white]
\tikzstyle{R}=[circle, inner sep=0, minimum size =6 pt, line width = 1pt, draw=black, fill=red, text= white]
\tikzstyle{c}=[circle, inner sep=0, minimum size =20 pt, line width = 1pt, draw=black, fill=white, text= black]
\tikzstyle{noeud_bleu}=[circle,inner sep=2, minimum size =3 pt, line width = 1pt, draw=black, fill=blue!60]
\tikzstyle{noeud_rouge}=[circle,inner sep=2, minimum size =3 pt, line width = 1pt, draw=black, fill=red!70]
\title{Incidence, a Scoring Positional Game on Graphs\footnote{This research was supported by the ANR project P-GASE (ANR-21-CE48-0001-01).}}

\author[$\dagger$]{Guillaume Bagan}
\author[$\dagger$]{Quentin Deschamps}
\author[$\dagger$]{Eric Duchêne}
\author[$\ddag$]{Bastien Durain}
\author[$\dagger$]{Brice Effantin}
\author[$*$]{Valentin Gledel}
\author[$\dagger$]{Nacim Oijid}
\author[$\dagger$]{Aline Parreau}

\affil[$\dagger$]{Univ. Lyon, Universit\'e Lyon 1, LIRIS UMR CNRS 5205, F-69621, Lyon, France.}
\affil[$\ddag$]{\'Ecole Normale Sup\'erieure de Lyon, 69364 Lyon Cedex 07, France.}
\affil[$*$]{Department of Mathematics and Mathematical Statistics, University of Ume\aa\,Sweden.}

\date{}

\begin{document}

\maketitle

\begin{abstract}

Positional games have been introduced by Hales and Jewett in 1963 and have been extensively investigated in the literature since then. These games are played on a hypergraph where two players  alternately select an unclaimed vertex of it. In the Maker-Breaker convention, if Maker manages to fully take a hyperedge, she wins, otherwise, Breaker is the winner. In the Maker-Maker convention, the first player to take a hyperedge wins, and if no one manages to do it, the game ends by a draw. In both cases, the game stops as soon as Maker has taken a hyperedge. By definition, this family of games does not handle scores and cannot represent games in which players want to maximize a quantity. 

\smallskip

In this work, we introduce scoring positional games, that consist in playing on a hypergraph until all the vertices are claimed, and by defining the score as the number of hyperedges a player has fully taken. We focus here on \inc, a scoring positional game played on a 2-uniform hypergraph, i.e. an undirected graph. In this game, two players alternately claim the vertices of a graph and score the number of edges for which they own both end vertices. In the Maker-Breaker version, Maker aims at maximizing the number of edges she owns, while Breaker aims at minimizing it. In the Maker-Maker version, both players try to take more edges than their opponent.

\smallskip

We first give some general results on scoring positional games such that their membership in Milnor's universe and some general bounds on the score. We prove that, surprisingly, computing the score in the Maker-Breaker version of \inc{} is {\sf PSPACE}-complete whereas in the Maker-Maker convention, the relative score can be obtained in polynomial time. In addition, for the Maker-Breaker convention, we give a formula for the score on paths by using some equivalences due to Milnor's universe. This result implies that the score on cycles can also be computed in polynomial time.
\end{abstract}

\section{Introduction}

\subsection{Positional games}

Positional games have been introduced by Hales and Jewett in 1963~\cite{Hales1963} and popularized by Erd{\H o}s and Selfridge in 1973 \cite{erdos1973}. Interest in them has increased due to the large number of games they can handle. 

In the standard definition of positional games, the board is a hypergraph on which two players alternately select an unclaimed vertex. In the Maker-Breaker convention, if Maker manages to claim all the vertices of a hyperedge, he wins, otherwise, Breaker is the winner. In the Maker-Maker convention, the first player, if any, who takes a hyperedge wins. If no player manages to claim all the vertices of a hyperedge, the game ends by a draw. 

Maker-Maker games are often considered as harder than Maker-Breaker games, since the objective of trying to fill a hyperedge and controlling at the same time that the opponent does not win, is hard to meet.

Positional games are finite perfect information two-players games. As such, there exists a winning strategy for one of the players or both players can insure a draw. The main issue is then to compute, for a given hypergraph, which player has a winning strategy. This problem has been proven to be {\sf PSPACE}-complete for both conventions, even if all the hyperedges have size at least 11 by Schaefer \cite{schaefer1978}. This result was recently improved to hypergraph with hyperedges of size at least 6 by Rahman and Watson \cite{rahman2021}.
On the other side, Galliot {\em et al.} proved that the winner can be computed in polynomial time on $3$-uniform hypergraphs \cite{galliot2022}.

In practice, positional games are studied in specific hypergraphs. Historically, they are almost always derived from hypergraphs built from a grid or a complete graph (see for example the reference books \cite{Beck2008,Hefetz2014}). 
More recently, some positional games played on hypergraphs derived from general graphs have been studied. For such games, Maker aims at building a structure in a given graph, and Breaker aims at preventing him to do so. The structure could be, for example, a copy of a graph $H$ ({\sc H-Game} \cite{Kronenberg2019}) or a dominating set ({\sc Maker-Breaker domination game} \cite{duchene2020}).

\subsection{Scoring games}

In parallel to the study of positional games, scoring games have been introduced in the 1950s by Milnor \cite{milnor1953} and Hanner \cite{hanner1959}. Their study was almost forgotten until the 2000s, when different formalisms for such games have been introduced by Ettinger~\cite{ettinger1996}, Stewart \cite{stewart2012}, or Larsson, Nowakowski and Santos~\cite{larsson2015}. The survey paper~\cite{larssonGONC} summarizes these different approaches.

In scoring games, two players, usually Left and Right, alternate moves with a score adjoined to the game. Each move of a player can modify this score, Left aims at maximizing the score at the end of the game, while Right tries to minimize it. Since scoring games are also finite perfect information games, if both players play optimally, the score at the end of the game is well-defined and only depends on who starts.

Despite the fact that scoring games were less studied, mainly due to the difficulty to build a general framework for them, particular scoring games on graphs have still been introduced recently. One can cite the game {\sc Influence} introduced by Duchêne {\em et al.} in 2021 \cite{duchene2021} which has been proven {\sf PSPACE}-complete in 2022 \cite{duchene2022}, or the largest connected subgraph game, introduced by Bensmail {\em et al.}, firstly as a scoring connection game \cite{bensmail2022}, and then as a Maker-Breaker connection game \cite{bensmail2022-2}. In \cite{larssonGONC}, there is a list of other particular scoring games on graphs that have been recently studied.

\subsection{Scoring positional games and outline of the paper}

In the current paper, we introduce a general scoring version of positional games. Left and Right alternately select vertices of a hypergraph until all the vertices are selected. Points are given when a hyperedge is fully selected by a player. In the Maker-Maker convention, both players get points and the score is the difference between the number of hyperedges taken by Left and  Right. In the Maker-Breaker version, the score is only the number of hyperedges taken by Left.

\paragraph{Outline of the paper.} In Section 2, after giving a formal definition of these games, we provide some general results on them. In particular, we prove that they belong to Milnor's universe and that determining the score is {\sf PSPACE}-complete in the two conventions. 
In the rest of the paper, we explore the game \inc{} that corresponds to the subcase of 2-uniform hypergraphs (or equivalently to graphs). In Section 3, we prove that, unlike for standard positional games, the Maker-Maker version of \inc{} is the easiest one since computing the score is linear in this case. Then we focus on the Maker-Breaker version of \inc. In Section 4, we give some general bounds on the score as well as some nice properties to deal with twin vertices. This allows us to calculate the exact value of the score for complete binary trees.
The next section shows that computing the score in Maker-Breaker convention is {\sf PSPACE}-complete but fixed-parameter-tractable when parameterized by the {\em neighbourhood diversity} of the graph (introduced in \cite{lampis2010}), which implies in particular that it is also FPT when parameterized by vertex cover. The last section is dedicated to the study of paths and cycles. We prove some equivalence relations between paths, which lead to a closed formula for paths and cycles. In particular, we can compute exactly the score for a path of length $n$, which is equal to $n/5+c$ where $c$ only depends on $n \bmod 5$.

\section{General results on scoring positional games}\label{scoring pos games}

\subsection{Definitions}\label{sec:defscore}

Scoring positional games are played on hypergraphs by two players, Left and Right, with the same rules as for standard positional games. The only difference lies in the winning convention. In a scoring positional game, the game ends when all vertices have been claimed. The score of a player is then defined as the number of hyperedges he manages to take. In the Maker-Maker convention, each player tries to maximize his score. In the Maker-Breaker convention, Maker (identified as Left) tries to maximize her score while Breaker (identified as Right) aims at minimizing the score of Maker. 

More formally, as for any scoring game, two scores are defined depending on which player starts. Let $H=(V,E)$ be a hypergraph. We define the score of $H$ as follows:  

\begin{itemize}
    \item in the Maker-Maker convention, $Ls(H)$ (resp. $Rs(H)$) as the difference between the scores of Left and Right when Left starts (resp. when Right starts) and both players play optimally.
    \item in the Maker-Breaker convention, $Ls(H)$ (resp. $Rs(H)$) as the score of Left when Left (resp. Right) starts and both players play optimally. 
\end{itemize}

It is well-known in scoring game theory that these notions exist and are well-defined (by considering the game tree of all the possible moves). 
Note that in the Maker-Maker convention, by symmetry of the roles of both players, we have that $Ls(H)=-Rs(H)$, so computing $Ls(H)$ will be of sufficient interest. In the Maker-Breaker convention, we have that $Ls(H)$ and $Rs(H)$ are nonnegative values by definition. 

\smallskip

In addition, it will be helpful to consider the scores obtained after some vertices have been claimed. A {\em position} of a scoring positional game is a triplet $P=(H,V_L,V_R)$ such that $V_L$ and $V_R$ are disjoint subsets of vertices. The set $V_L$ corresponds to the vertices claimed by Left whereas $V_R$ correspond to the vertices claimed by Right. The set of remaining vertices will be generally denoted by $V_F$. We have $V_F=V\setminus (V_L\cup V_R)$. For both conventions, we will denote by $Ls(P)$ (resp. $Rs(P)$) the score of $H$ if Left has already claimed the vertices of $V_L$, and Right the vertices of $V_R$, when Left (resp. Right) starts.
When $V_F\neq \emptyset$, the scores at a position $P$ can be recursively defined as follows:
\begin{align*}
    &Ls(P) = \underset{x\in V_F}{\max} Rs(H,V_L\cup\{x\},V_R)\\
    &Rs(P) = \underset{x\in V_F}{\min} Ls(H,V_L,V_R\cup\{x\}).\\
\end{align*}
When $V_F=\emptyset$, the score depends on the convention. In Maker-Maker convention,  $$Ls(P)=Rs(P)=|\{e\in E | e\subseteq V_L\}|-|\{e\in E |e\subseteq V_R\}|$$ whereas in Maker-Breaker convention, we have $$Ls(P)=Rs(P)=|\{e\in E | e\subseteq V_L\}|.$$

\smallskip

In the literature, there are few games that can be seen as scoring positional games. The famous Dots and Boxes games \cite{berlekamp2000}, that has recently be proven {\sf PSPACE}-complete by Buchin {\em et al.} \cite{buchin2021}, could be an example, with the additional constraint that a player is forced to move again each time he gets points. By removing this constraint, we get a pure example of the above definition (in the Maker-Maker convention), and the game is known as Picarête \cite{blanc2006}. More recently, the Constructor-Blocker game introduced by Patkos {\em et al.} \cite{patkos2022} in 2022, in which Constructor aims at maximizing the number of copies of a graph $H$ with a forbidden graph $F$, can be seen as a scoring positional game when $F$ is empty.

\paragraph{Incidence} In most of this paper, we will mainly focus on an example of scoring positional game that is called \inc. It corresponds to the game played on a hypergraph where all hyperedges are of size two. In others terms, this game can be defined as follows on a simple graph $G=(V,E)$.
Alternately, two players claim an unclaimed vertex of $V$. When all the vertices have been taken, the score of a player is defined as the number of edges in the subgraph of $G$ induced by the vertices he claimed. 

Hence, in both conventions, Left (that is always Maker) aims at collecting points by claiming the two extremities of an edge. The main difference concerns the role of Right, that aims at touching the maximum number of edges (hence prohibiting a maximum number of points for Left)  in the Maker-Breaker convention. See Figure \ref{fig:score_conventions} for an example of computations of the score at the end of a game.

\begin{figure}[ht]
\centering
    \begin{tikzpicture}
        \node[noeud_rouge,label=below:$R$] (a) at (0,0){};
        \node[noeud_bleu,label=below:$L$] (b) at (1.5,0){};
        \node[noeud_rouge,label=above:$R$] (c) at (0,1.5){};
        \node[noeud_rouge,label=above:$R$] (d) at (1.5,1){};
        \node[noeud_bleu,label=below:$L$] (e) at (2.75,0.5){};
        \node[noeud_bleu,label=above:$L$] (f) at (2.4,2){};
        \node[noeud_bleu,label=below:$L$] (g) at (3.75,1.5){};
        
        \draw (a) -- (b);
        \draw (a) -- (c);
        \draw (b) -- (d);
        \draw (b) -- (e);
        \draw (c) -- (d);
        \draw (c) -- (f);
        \draw (d) -- (f);
        \draw (e) -- (f);
        \draw (e) -- (g);
        \draw (f) -- (g);
    \end{tikzpicture}
\caption{\label{fig:score_conventions}An endgame of \inc. In Maker-Maker convention the score of the position is 2 while it is 4 in Maker-Breaker convention.}
\end{figure}

\subsection{Milnor's universe}\label{sec:milnor}

In 1953  \cite{milnor1953}, Milnor introduced a universe of scoring games having nice properties. This universe is the one of {\it dicotic nonzugzwang} games: 
\begin{itemize}
\item a game is {\it dicotic} if at any moment of the game, if a player can play, the other player can also play.
\item a game is {\it nonzugzwang} if at any moment of the game, both players have no interest in skipping their turn. In the context of scoring positional games, it means that for a hypergraph $H$, we have $Ls(H,V_L, V_R) \ge Rs(H,V_L, V_R)$ for any sets of vertices $V_L, V_R$ claimed by Left and Right during the game.
\end{itemize}

Being in Milnor's universe induces a couple of useful results concerning the sum operator and the equivalence of games. The {\em disjunctive sum operator} $+$ applied to scoring (positional) games $G_1$ and $G_2$ defines the game $G_1+G_2$ as the game in which a move consists in moving either in $G_1$ or in $G_2$. The game ends when the moves in both components of the sum are exhausted. See \cite{duchene2022} for the formal definition.
Note that the sum of two scoring positional games, with the same convention, is still a scoring positional game with hypergraph the disjoint union of the two hypergraphs. 
As game sums appear in many games when playing, one could expect to simplify them by replacing large games by smaller ones. This leads to the notion of equivalence of games:

\begin{definition}[Milnor \cite{milnor1953}]
Two scoring games $G_1$ and $G_2$ are \em equivalent \em (write $G_1 \equiv G_2$) if for any game $G$, we have $Ls(G + G_1) = Ls(G + G_2)$ and $Rs(G+G_1) = Rs(G+G_2)$.
\label{milnor equiv}
\end{definition}

In other terms, one can always exchange $G_1$ and $G_2$ in any sum of games if they are equivalent. In particular, games that are equivalent to the empty game can be removed from any sum of games.\\

Games belonging to Milnor's universe form an Abelian group with the sum operator\cite{milnor1953}. In particular, this implies that every game $G$ in Milnor's universe admits an inverse, i.e. a game $G'$ such that $G+G' \equiv 0$ (where $0$ is the empty game). More precisely, this inverse corresponds to the {\em negative} of $G$, i.e. the game where the roles of Left and Right are exchanged, together with their scores.\\ 

Moreover, proving equivalence in Milnor's universe is greatly simplified, thanks to the next lemma.

\begin{lemma}[Milnor \cite{milnor1953}]\label{lemma equiv}
For any games $G$ and $H$ that are dicotic nonzugzwang, we have:
$Ls(G-H)=Rs(G-H)=0$ if and only if $G$ and $H$ are equivalent. 
\end{lemma}

In addition, sums of games in Milnor's universe can be bounded as follows:
\begin{lemma}[Milnor \cite{milnor1953}]\label{lemma encadrement}
Let $G$ and $H$ be two dicotic nonzugzwang games, we have 
$$
Rs(G) + Rs(H) \leq  Rs(G+H) \leq Ls(G)+Rs(H)\leq Ls(G+H)\leq Ls(G)+Ls(H).
$$
\end{lemma}

In what follows, we will show that scoring positional games belong to Milnor's universe. Yet, the negative of a game cannot be defined in the Maker-Breaker convention, as the scores of Maker and Breaker can not be interchanged naturally, by asymmetry of the definition of the score. Therefore, we have decided to embed scoring positional games in a more general family that will be called {\sc partisan scoring positional games}. The term partisan is derived from standard combinatorial games~\cite{berlekamp2000}, meaning that Left and Right may have different moves (and also different ways of scoring points). \\

A partisan scoring positional game is played on a hypergraph $H$ whose hyperedges are either colored blue, red or green. The two players, Left and Right, alternatively claim vertices of $H$. The score of Left corresponds to the blue and green hyperedges she claimed, whereas the score of Right corresponds to the red and green ones. As previously, the score of the game ($Ls(H)$ and $Rs(H)$, depending on who starts) is the difference between the score of Left and Right.\\

Partisan scoring positional games include both Maker-Maker and Maker-Breaker scoring positional games. Even more, the convention can be omitted, as it is deduced by the colors of the hypergraph. Indeed, if all the hyperedges are green, it means that both players can win any hyperedge, which corresponds to the Maker-Maker version. If all the hyperedges are blue, it corresponds to the Maker-Breaker convention, as only Left can get points.
According to this definition, the negative of a partisan scoring positional game is well-defined, as it suffices to exchange the colors blue and red in the hyperedges, as well as the vertices already chosen by Left and Right (if any).

We will now give several general results about partisan scoring positional games. By inclusion, these results will also concern scoring positional games. First, we will prove that they belong to Milnor's universe and thus satisfy Lemma~\ref{lemma equiv}.

\begin{lemma}
Partisan scoring positional games belong to Milnor's universe.
\end{lemma}

\begin{proof}
Let $H = (V,E)$ be a hypergraph with hyperedges colored blue, red and green, and $V_L, V_R \subset V$ be vertices already claimed by Left and Right respectively such that $V_L \cap V_R = \emptyset$. \\

\noindent {\em A partisan scoring positional game is dicotic:} if $V_L \cup V_R = V$, then no moves are available, neither for Left nor for Right. Otherwise, let $v \in V \setminus \{V_L \cup V_R\}$ . Both Left and Right are allowed to play $v$ as it is an unclaimed vertex. Therefore, the game is dicotic.\\

\noindent {\em A partisan scoring positional game is nonzugzwang:} We need to prove that $Ls(H, V_L, V_R) \ge Rs(H, V_L, V_R)$. Let $k = Rs(H, V_L, V_R)$ with $V_L, V_R$ vertices already claimed in $H$ by Left and Right respectively. If $V_L\cup V_R = V$, we have $Ls(H, V_L, V_R) = Rs(H, V_L, V_R)= k$ as there is no move available in $H$. Otherwise, let $\mathcal{S}$ be an optimal strategy for Left when Right starts. We define a strategy $\mathcal{S}'$ for Left when she starts as follows:

\begin{itemize}
    \item Left considers an arbitrary unclaimed vertex $v_0$ of the graph, and plays the vertex she would have played in $\mathcal{S}$ if Right plays $v_0$.
    \item Whenever, Right plays a vertex $w$ in $V\setminus\{v_0\}$, she plays the vertex she would have played in $\mathcal{S}$ if Right has played $w$ in $\mathcal{S}$ after having played $v_0$ on first move.
    \item If Right plays $v_0$, she considers an arbitrary unclaimed vertex $v_1$ in the graph, and continues this strategy by supposing that Right has played $v_1$ instead of $v_0$. More generally, when Right claims the vertex $v_{\ell}$, she considers an unclaimed vertex $v_{\ell+1}$ and considers that Right has claimed $v_{\ell+1}$ instead. 
    \item At the end, if she needs to consider that Right has played a vertex $v_{\ell}$ and no other vertex is available, she plays $v_{\ell}$.
\end{itemize}

Following this strategy, all the vertices Left would have played in $\mathcal{S}$ if Right has played the vertices $v_i$s she has considered, have been played in $\mathcal{S'}$ by Left. Similarly, the vertices that Right have played in $\mathcal{S}'$ are a subset of the one he would have played in $\mathcal{S}$. Therefore, as $\mathcal{S}$ was an optimal strategy in $H$ when Right starts, this strategy ensures that Left scores at least $k = Rs(H, V_L, V_R)$. Finally, we have $Ls(H, V_L, V_R) \ge k = Rs(H, V_L, V_R)$, and the game is nonzugzwang.

As the game is nonzugzwang and dicotic, it belongs to Milnor's universe.
\end{proof}

As a consequence, this result applies also to scoring positional games and, in particular, the game \inc. We will use this result in Section~\ref{section path} to solve \inc{} on paths.

\subsection{Algorithmic complexity}

We now prove that deciding if the Left scores of a scoring positional game is {\sf PSPACE}-complete in both conventions. This result is a direct consequence of the {\sf PSPACE}-complexity of (non-scoring) positional games.
\medskip

\noindent \decisionpb{\MBPG}
{A hypergraph $H$, $P\in \{\text{Maker},\text{Breaker}$\}.}
{True if Maker wins the Maker-Breaker positional game played on $H$ with first player $P$.}{1}\\[0.5em]

\MBPG{} has been proved to be {\sf PSPACE}-complete by Schaeffer \cite{schaefer1978} for 11-uniform hypergraphs (all the hyperedges have size 11). This result was recently improved to 6-uniform hypergraphs by Rahman and Watson \cite{rahman2021}.

\begin{theorem}[\cite{rahman2021}]\label{thm:pspaceMBPG}
\MBPG{} is {\sf PSPACE}-complete even restricted to $6$-uniform hypergraphs. 
\end{theorem}

\MBPG{} can easily be reduced to the two following problems on scoring positional games.
\medskip

\noindent \decisionpb{\MBSPG}{A hypergraph $H$, an integer $k$, a first player $P\in \{Left,Right\}$.}
{True if the $P$ score in the scoring positional game played on $H$ with Maker-Breaker convention is at least $k$, false otherwise.}{1}

\medskip

\noindent \decisionpb{\MMSPG}{A hypergraph $H$, an integer $k$.}
{True if the Left score in the scoring positional game played on $H$ with Maker-Maker convention is at least $k$, false otherwise.}{1}

\begin{corollary}\label{cor:generalcomplexity}
\MBSPG{} is {\sf PSPACE}-complete even restricted to $6$-uniform hypergraphs, $P=Left$ and $k=1$.

\MMSPG{} is {\sf PSPACE}-complete even restricted to $7$-uniform hypergraphs and $k=1$.
\end{corollary}

\begin{proof}
Since both games are played in $|V(H)|$ turns, they belong to {\sf PSPACE} according to Section~6.1 in \cite{Hearn2009}.

Let $H$ be a $6$-uniform hypergraph and assume Left is the first player. We have $Ls(H) \ge 1$ in the Maker-Breaker convention if and only if Maker wins the Maker-Breaker positional game (without score) played on $H$ with Maker as first player. Thus, by Theorem \ref{thm:pspaceMBPG}, \MBSPG{} is {\sf PSPACE}-complete even restricted to $6$-uniform hypergraphs, $k=1$ and $P=Left$.

Consider now $H'$ the $7$-uniform hypergraph obtained from $H$ by adding a universal vertex $v_0$: each hyperedge of $H$ is extended to contain $v_0$. There exists an optimal strategy in the Maker-Maker convention that starts by claiming $v_0$.
Then the other player cannot score any point. Then, we have $Ls(H') \ge 1$ if and only if Maker wins playing second in the Maker-Breaker positional game (without score) played on $H$. Thus, by Theorem \ref{thm:pspaceMBPG}, \MMSPG{} is {\sf PSPACE}-complete even restricted to $7$-uniform hypergraphs and $k=1$.
\end{proof}

We will complete the results of Corollary \ref{cor:generalcomplexity} in next sections by proving that \MBSPG{} is still {\sf PSPACE}-complete for $2$-uniform hypergraphs (Theorem \ref{thm:pspaceincidence}). This will imply that \MMSPG{} is {\sf PSPACE}-complete for $3$-uniform hypergraphs. To complete the picture, we will give a linear algorithm to solve \MMSPG{} in $2$-uniform hypergraphs (Theorem \ref{thm:polyMM}).

\subsection{Bounds in Maker-Maker convention}

In this subsection, we give an easy bound on the score in Maker-Maker convention, using the maximal degree of the hypergraph. Let $H$ be a hypergraph. The {\em degree} of a vertex $v$ of $H$ is the number of hyperedges containing $v$. We denote by $\Delta(H)$ the maximal degree of $H$.

\begin{lemma}\label{positive}
Let $H$ be a hypergraph. In the Maker-Maker scoring positional game on $H$, we have $- \Delta(H) \le Rs(H) \le 0 \le Ls(H) \le \Delta(H)$.
\end{lemma}

\begin{proof}
As noticed in Section \ref{sec:defscore}, we have $Ls(H)=-Rs(H)$ in the Maker-Maker convention since players have symmetric roles.
Since the game is nonzugzwang, we also have $Ls(H) \ge Rs(H)$ which implies that  $Rs(H) \le 0 \le Ls(H)$.

To prove the upper bound with $\Delta(H)$, we just need to prove that $Ls(H) \le \Delta(H)$. Let $v_0$ be the first vertex played in an optimal strategy. Consider the hypergraph $H'$ obtained from $H$ by removing $v_0$ and all the hyperedges containing it. If the second player applies the optimal strategy for $H'$ during the rest of the game, he will score at least $Rs(H') \le 0$ on it and the final score will be at most $|\{ e | v_0 \in e\}| + Rs(H')$. Thus, we have $Ls(H)\le deg(v_0)+Rs(H') \le \Delta(H)$.
\end{proof}

We do not think that the upper bound in Lemma \ref{positive} is tight if the hypergraph is simple (i.e. there are no two hyperedges that contain exactly the same vertices). Actually, the best example we know in this case is a hypergraph $H$ having a universal vertex $x$, a hyperedge with $x$ alone and $\Delta-1$ hyperedges of size 2 containing $x$ and another unique vertex, see Figure~\ref{fig:boundDelta}. For this hypergraph, $Ls(H)=\lfloor \frac{\Delta(H)+1}{2}\rfloor$. 
Besides, we will prove that for $2$-uniform hypergraphs (i.e. graphs),  the score is at most $\Delta(H)/2$ (see Corollary \ref{cor:ubMMdelta}. We believe that this bound remains true in any hypergraph:

\begin{figure}
    \centering
        \begin{tikzpicture}
\node[v] (A) at (0,0) {};
 \node[v] (B) at (60:1.5) {};
 \node[v] (C) at (0:1.5) {};
 \node[v] (D) at (120:1.5) {};
 \node[v] (E) at (180:1.5) {};

\draw (180:0.75) ellipse (1cm and 0.3cm);
\draw (0:0.75) ellipse (1cm and 0.3cm);
\draw[rotate around={60:(60:0.75)}] (60:0.75) ellipse (1cm and 0.3cm);
\draw[rotate around={120:(120:0.75)}] (120:0.75) ellipse (1cm and 0.3cm);

 \draw (0,0) circle (0.4);
    \end{tikzpicture}

    \caption{A hypergraph satisfying $Ls(H) = \lfloor \dfrac{\Delta(H)+1}{2}\rfloor$ in Maker-Maker convention}
    \label{fig:boundDelta}
\end{figure}
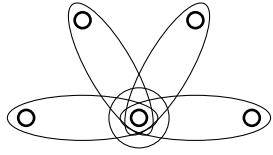

\begin{conjecture}
Let $H$ be a simple hypergraph. In the Maker-Maker scoring positional game on $H$, we have $Ls(H) \le \dfrac{\Delta(H)+1}{2}$.
\end{conjecture}

\subsection{Bounds in Maker-Breaker}

In Maker-Breaker convention, the bound from Lemma \ref{positive} is not valid anymore. Indeed, the score can actually be linear with the number of vertices of the hypergraph, even if the maximal degree is constant. Next, we derive a general tight bound, based on the same principle used to prove the Erdös-Selfridge criterion \cite{erdos1973}. Some tight examples will be given in Section \ref{general result MB} for $2$-uniform hypergraphs (see Corollary~\ref{cor:ESgraph}).

\begin{theorem}[Erd{\H o}s, Selfridge, 1973 \cite{erdos1973}]
Let $H = (V,E)$ be a hypergraph. If $\underset{e \in E}{\sum} 2^{-|e|} < 1$, then Breaker wins on $H$ when he starts. If $\underset{e \in E}{\sum} 2^{-|e|} < \frac{1}{2}$, then Breaker wins on $H$ when Maker starts. 
\end{theorem}

The main idea to prove this theorem is that if the hyperedges are large enough, Breaker will have the time to play in all of them before Maker can fill one. A similar idea can be introduced when dealing with scores by computing how many hyperedges Breaker can touch. The strategy used relies on a greedy strategy by introducing a potential function, as it was done by Erd{\H o}s and Selfridge. Let $H$ be a hypergraph. We denote by $\ell(H)$ the maximum number of hyperedges that contain a fixed pair of vertices. More formally, $\ell(H)= \underset{x,y \in V^2}{\max}|\{e \in E | x,y\in e\}|$.

\begin{theorem}\label{ES like bound}
Let $H = (V,E)$ be a hypergraph. In the Maker-Breaker convention, we have $Ls(H) \ge \underset{e \in E}{\sum} 2^{-|e|} - \frac{n\ell(H)}{8}$, and $Rs(H) \le \underset{e \in E}{\sum} 2^{-|e|}$.
\end{theorem}

\begin{proof}
Let $(H,V_L,V_R)$ be any position of a Maker-Breaker scoring positional game.
We introduce the potential function: $$P(H,V_L,V_R) = \underset{e \in E, e\cap V_R=\emptyset}{\sum} 2^{-|e\setminus V_L|}.$$
In this function, only hyperedges not played by Right are considered, and we only count the number of free vertices in the edge. Note that at the beginning of the game, $P(H,\emptyset,\emptyset)=\underset{e \in E}{\sum} 2^{-|e|}$. At the end of the game, $V=V_L\cup V_R$ and $P(H,V_L,V_R)= |\{e\in E | e \cap V_R=\emptyset\}|$ is the final score. Furthermore, when a vertex $v$ is played by Maker (respectively Breaker), the potential is increasing (resp. decreasing) by the quantity $$\delta_P(H,V_L,V_R,v)=\sum_{e|e\cap V_R=\emptyset,v \in e}2^{-|e\setminus V_L|}.$$

Let $\mathcal{S}$ be a strategy for Maker consisting in maximizing $P$ at each move, i.e. Maker chooses the vertex $v$ that maximizes $\delta_P(H,V_L,V_R,v)$. We prove that this strategy provides the desired bound. Suppose first that Maker starts. Suppose $V_L$ and $V_R$ have already been played by Maker and Breaker respectively. Let $v_L$ the vertex played by Maker according to $\mathcal{S}$ and $v_R$ the vertex played by Breaker after this move. As Maker has played $v_L$ and not $v_R$, we have, before $v_L$ was played, $\delta_P(H,V_L,V_R,v_L)\geq \delta_P(H,V_L,V_R,v_R)$.

However, $\delta_P(H,V_L\cup\{v_L\},V_R,v_R)$ might be larger than $\delta_P(H,V_L,V_R,v_R)$ after $v_L$ was played if there exist some hyperedges that contain both $v_L$ and $v_R$. We actually have:

\begin{align*}
\delta_P(H,V_L\cup\{v_L\},V_R,v_R)&=\delta_P(H,V_L,V_R,v_R)+\sum_{e\cap V_R=\emptyset, v_L,v_R\in e}2^{-|e\setminus V_L|}\\
&\leq \delta_P(H,V_L,V_R,v_R) + \frac{\ell(H)}{4}.
\end{align*}

Last inequality comes from the fact that $e\setminus V_L$ must contain $v_L$ and $v_R$ and thus has size at least 2. 
Therefore, we have 
\begin{align*}
P(H, V_L \cup \{v_L\}, V_R \cup \{v_R\})&= P(H,V_L,V_R) +\delta_P(H,V_L,V_R,v_L)-\delta_P(H,V_L\cup\{v_L\}, V_R,v_R)\\
&\geq P(H,V_L,V_R)- \frac{\ell(H)}{4}.
\end{align*}
As there is $n$ moves in the game by applying this step $\frac{n}{2}$ times for each pair of moves (recall that we consider here that Maker starts), we have at the end of the game $Ls(H)\geq P(H, V_L, V_R) \ge P(H, \emptyset, \emptyset) - \frac{n}{2}\times\frac{\ell(H)}{4}$, as required.

Suppose now that Breaker starts and considers this strategy for him (i.e. choosing the vertex $v$ that maximizes $\delta_P(H,V_L,V_R,v)$). Suppose $V_L$ and $V_R$ have already been played by Maker and Breaker respectively. Let $v_R$ be the vertex played by Breaker according to $\mathcal{S}$ and let $v_L$ be the vertex answered by Maker. We have 
$\delta_P(H,V_L,V_R,v_R) \geq \delta_P(H,V_L,V_R,v_L)$.
 Note that here, $\delta_P(H,V_L,V_R\cup\{v_R\},v_L)$ cannot increase after the move of Right, as it does not change the size of the hyperedges (it can only decrease if some edges containing $v_L$ also contains $v_R$). Therefore, after these two moves, we obtain $ P(G, V_L\cup \{v_L\}, V_R \cup \{v_R\}) \le P(G,V_L,V_R)$. By applying this result from $V_L = V_R = \emptyset$ to the end of the game, we obtain $P(H,V_L,V_R) \le P(H, \emptyset, \emptyset)$ for any sets $V_L$ and $V_R$ obtained after Right applies $\mathcal{S}$. In particular, when the game ends, this strategy ensures that $Rs(H) \le P(H, \emptyset, \emptyset) = \underset{e \in H}{\sum} 2^{-|e|}$.
\end{proof}

From now on and until the end of the paper, we will focus on the game \inc, i.e. the scoring positional game played on a graph.

\section{{\sc Maker-Maker Incidence} is polynomial}

In this section, we provide a linear time algorithm  to compute the score of \mminc. A natural idea, while playing \inc, is that high degree vertices are interesting to play first, as they enable to score many points with their multiple adjacent edges. Therefore, a simple strategy for both players would be to play greedily by always picking an available vertex of highest degree. We here prove that this strategy is optimal.

Later in Section~\ref{pspace}, we will prove that {\mbinc} is {\sf PSPACE}-complete, which induces that {\MMSPG} is {\sf PSPACE}-complete on $3$-uniform hypergraphs.

\begin{theorem}\label{thm:polyMM}
Let $G$ be a graph with $n$ vertices. Let $d_1\geq ... \geq d_n$ be the degree of the vertices in decreasing order. For the game {\mminc} played on $G$, we have
$$
Ls(G)=\frac{1}{2} \left (\sum_{i \text{ odd }}  d_{i} - \sum_{i \text{ even }}  d_{i}\right ).
$$

In particular, the score can be computed in linear time.
\end{theorem}

\begin{proof}
Let $G = (V,E)$ be a graph. Denote by $v_1, \dots, v_{n}$ the vertices of $G$ of degree $d_1$, \dots, $d_{n}$ respectively, and arranged such that $d_1 \ge d_2 \ge \dots \ge d_{n}$.  
Denote by $s = \frac{1}{2}  (\underset{\tiny i \text{ odd}}{\sum} d_{i} - \underset{i \text{ even}}{\sum}  d_{i})$. We will prove that $Ls(G) = s $. Before proving the value of the score, we prove the following claim:\\

\begin{claim}
 Denote by $V_L$ the vertices claimed by Left, and by $V_R$ the vertices claimed by Right at the end of a game played on $G$. The score obtained is $\frac{1}{2} ( \underset{v_l \in V_L }{\sum} d_l - \underset{v_r \in V_R }{\sum} d_r )$.
\end{claim}

\begin{claimproof}
Denote by $e_L$ (resp. $e_R$) the number of edges where both endpoints were claimed by Left (resp. Right) and by $e_0$ the number of edges which have one extremity claimed by each player.

By definition, the score is $e_L - e_R$. Now, by a double counting argument, we have $\underset{v_l \in V_L }{\sum} d_l = 2 e_L + e_0 $, and $\underset{v_r \in V_R  }{\sum} d_r = 2 e_R + e_0 $. Therefore, the score of the game is $e_L - e_R = \frac{1}{2} ( \underset{v_f \in V_L }{\sum} d_l - \underset{v_r \in V_R }{\sum} d_r )$.
\end{claimproof}

Now we provide a strategy for Left that proves that $Ls(G) \ge s$. The same argument works for Right and leads to $Ls(G) \le s$. Consider that Left claims at each turn the free vertex of highest degree. During her first turn, she claims a vertex of degree $d_1$, during the second turn, she claims either a vertex of degree $d_2$ or $d_3$, both having a value of at least $d_3$, \dots, during here $k$-th turn, she will claim a vertex of degree $d_{k}, d_{k+1}, \dots$ or $d_{2k-1}$, each of them have a value of at least $d_{2k-1}$. In the end, she will have played $\left \lceil \frac{n}{2}\right \rceil$ vertices, and the $k$-th of them will be of degree at least $d_{2k-1}$. Reciprocally, the highest degree played by Right has value at most $d_2$, the second highest has value at most $d_4$ and so on. Therefore, by using the result of the claim, the score obtained by this strategy is at least $s$.

The above score can be computed in linear time because it does not require to sort the list of the vertices, but only to know the number of vertices of any degree, which is bounded by $n-1$. 
\end{proof}

\begin{corollary}
Let $n \in \mathbb{N}$. Denote by $P_n$ the path of order $n$. In \mminc, we have $Ls(P_n) = -Rs(P_n) = 0$ if $n$ is even and $Ls(P_n) = -Rs(P_n) = 1$ if $n$ is odd. 
\end{corollary}

\begin{proof}
$P_n$ has exactly $n-2$ vertices of degree $2$ and two vertices of degree $1$. Therefore, if $n$ is even, an optimal strategy gives $\frac{n}{2}-1$ vertices of degree two and one vertex of degree one to each player, which provides a draw.
If $n$ is odd, Left has one more vertex of degree $2$ to play, and her score is then $1$.
\end{proof}

\begin{corollary}\label{cor:ubMMdelta}
Let $G$ be a graph of maximal degree $\Delta$. In \mminc{}, we have $Ls(G) \le \frac{\Delta}{2}$.
\end{corollary}

\begin{proof}
Let $G$ be a graph of maximal degree $\Delta$. Up to add an isolated vertex, suppose it has an even number of vertices. Denote by $d_1, d_2, \dots, d_{2n}$ its degrees written in decreasing order. 
We have $Ls(G)=\frac{1}{2} \underset{i = 1}{\overset{n}{\sum}} (d_{2i-1} - d_{2i}) = \frac{\Delta}{2} - \underset{i = 1}{\overset{n}{\sum}} (d_{2i} - d_{2i+1}) $, by setting $d_{2n+1} = 0$. For any $1 \le i \le n$, we have $d_{2i} \ge d_{2i+1}$. Hence, each term of the sum is nonnegative, and finally, we have $Ls(G) \le \frac{\Delta}{2}$.
\end{proof}

\section{General results on  {\sc Maker-Breaker Incidence}}\label{general result MB}

In the rest of the paper, we focus on the Maker-Breaker version of \inc. Contrary to the Maker-Maker version of this game, a greedy strategy is not always optimal. Thus, studying this game is much more challenging. In this section, we give some general results on this version.
We start with a direct application of the bound given for general scoring positional games in Theorem \ref{ES like bound}.

\begin{corollary}\label{cor:ESgraph}
Let $G$ be a graph with $n$ vertices and $m$ edges. In the \mbinc{} game, $Ls(G) \ge \frac{m}{4} - \frac{n}{8}$, and $Rs(G) \le \frac{m}{4}$.

These bounds are tight.
\end{corollary}

\begin{proof}
This is a direct application of Theorem~\ref{ES like bound}. Since the hypergraph is $2$-uniform and simple, for each pair of vertices, there is at most one edge containing the two vertices. Thus we have $\ell(G)=1$. Furthermore, each edge has size 2, thus $\underset{e \in G}{\sum} 2^{-|e|} = \frac{m}{4}$.

For tightness, consider first a graph $G$ that is a complete graph of order $8k$, with $k\in \mathbb N$. The lower bound gives $Ls(G)\geq \dfrac{{8k \choose 2}}{4}  - k={4k \choose 2}$. 
By playing randomly, Left takes $4k$ vertices and each pair of vertices scores one point. Thus $Ls(G)={4k \choose 2}$

Consider the graph $H$ made by a disjoint union of $2k$ paths on three vertices.
Left playing second can take $k$ central vertices and one leaf for each central vertex  he has taken. This strategy gives at most $k$ points to Left which is equal to the upper bound $\frac{m}{4}$ given in the statement.
\end{proof}

While playing \inc, some moves are equivalent: playing one or the other will not change the final score. This is in particular the case when two vertices have the same neighbourhood (up to the vertices already played). An interesting fact in this case is that, in Maker-Breaker convention, we can assume that each player will take exactly one of the two vertices.
More formally, let $G = (V,E)$ be a graph and $P=(G,V_L, V_R)$ some position of the game on $G$. Let $v_1, v_2$ be two free vertices. Vertices $v_1, v_2$ are said to be {\em equivalent} in $P$ if and only if we have $N(v_1) \cap V_F \setminus \{v_2\} = N(v_2) \cap V_F \setminus \{v_1\}$ and $|N(v_1) \cap V_L| = |N(v_2) \cap V_L|$. Note that the first equality is a set equality, while the second one only is on cardinals.

\begin{lemma}\label{twin vertices}
Let $G = (V,E)$ be a graph and let $P=(G,V_L, V_R)$ be a position of the game. Let $v_1, v_2$ be equivalent vertices in $P$. In \mbinc{}, we have $Ls(P) = Ls(G, V_L \cup \{v_1\}, V_R \cup \{v_2\})$ and $Rs(P) = Rs(G, V_L \cup \{v_1\}, V_R \cup \{v_2\})$.
\end{lemma}

\begin{proof}
We prove both results by induction on $|V_F| = |V \setminus (V_L \cup V_R)|$, the number of free vertices.
The result is clear if there are only two free vertices $v_1$ and $v_2$ as each player will claim one of them, and they will have the same number of neighbors in $V_L$ at the end.
Let $P = (G, V_L, V_R)$ be a position with $|V_F| \ge 3$, and let $v_1, v_2 \in V_F $ be equivalent vertices in $P$. 

\vspace{0.2cm}

We first prove that $Ls(P) = Ls(G, V_L \cup \{v_1\}, V_R \cup \{v_2\})$.
Let $x$ be an optimal move for Left in $P$. If $x \in \{v_1, v_2\}$, we have $Ls(P) = Rs(G, V_L \cup \{v_1\}, V_R)$. Indeed, exchanging the roles of $v_1$ and $v_2$ is possible since they will score exactly the same number of points at the end. 
Using the recursive definition of the scores we have, $Rs(G,V_L\cup \{v_1\},V_R) \le Ls(G, V_L \cup\{v_1\}, V_R \cup \{v_2\})$.
Otherwise, we have $Ls(P) = Rs(G, V_L \cup \{x\}, V_R)$.  Vertices $v_1$ and  $v_2$ are still equivalent in $(G, V_L \cup \{x\}, V_R)$. By induction, $Rs(G, V_L \cup \{x\}, V_R)=Rs(G, V_L \cup \{v_1, x\}, V_R \cup \{v_2\})$. According to the recursive definition of the score, $Ls(G, V_L \cup \{v_1\}, V_R \cup \{v_2\})\geq Rs(G, V_L \cup \{v_1, x\}, V_R \cup \{v_2\})$.
Finally, in both cases, $Ls(P)\leq Ls(G, V_L \cup \{v_1\}, V_R \cup \{v_2\})$.

We now prove the other inequality. Let $x$ be an optimal move for Left in $(G, V_L \cup \{v_1\}, V_R \cup \{v_2\})$. We have $Ls(G, V_L \cup \{v_1\}, V_R \cup \{v_2\}) = Rs(G, V_L \cup \{v_1, x\}, V_R \cup \{v_2\})$.
By induction, since $v_1$ and $v_2$ are still equivalent in $(G,V_L \cup \{x\}, V_R)$, we have $Rs(G,V_L \cup \{v_1,x\}, V_R\cup\{v_2\}) = Rs(G, V_L \cup \{x\}, V_R)$. Using the recursive definition of the score,  $Ls(P)\geq Rs(G, V_L \cup \{x\}, V_R)$, which leads to $Ls(P)\geq Ls(G,V_L\cup\{v_1\},V_R\cup\{v_2\})$. Finally, we have proved $Ls(P) = Ls(G, V_L \cup \{v_1\}, V_R \cup \{v_2\})$.

\vspace{0.2cm}

We now turn to the proof of $Rs(P) = Rs(G, V_L \cup \{v_1\}, V_R \cup \{v_2\})$. Let $x$ be an optimal move for Right in $P$. If $x \in \{v_1, v_2\}$, we have $Rs(P) = Ls(G, V_L, V_R\cup \{v_2\})$. Indeed, exchanging the roles of $v_1$ and $v_2$ is possible since they will score exactly the same number of points at the end. 
Using the recursive definition of the scores, we have $Ls(G,V_L,V_R\cup \{v_2\}) \ge Rs(G, V_L \cup\{v_1\}, V_R \cup \{v_2\})$.
Otherwise, we have $Rs(P) = Ls(G, V_L , V_R\cup \{x\})$.  Vertices $v_1$ and  $v_2$ are still equivalent in $(G, V_L, V_R\cup \{x\})$. By induction, $Ls(G, V_L, V_R\cup \{x\})=Ls(G, V_L \cup \{v_1\}, V_R \cup \{v_2, x\})$. According to the recursive definition of the score, $Rs(G, V_L \cup \{v_1\}, V_R \cup \{v_2\})\leq Ls(G, V_L \cup \{v_1\}, V_R \cup \{v_2, x\})$.
Finally, in both cases, $Rs(P)\geq Rs(G, V_L \cup \{v_1\}, V_R \cup \{v_2\})$.

We now prove the other inequality. Let $x$ be an optimal move for Right in $(G, V_L \cup \{v_1\}, V_R \cup \{v_2\})$. We have $Rs(G, V_L \cup \{v_1\}, V_R \cup \{v_2\}) = Ls(G, V_L \cup \{v_1\}, V_R \cup \{v_2, x\})$.
By induction, since $v_1$ and $v_2$ are still equivalent in $(G,V_L, V_R \cup \{x\})$, we have $Rs(G,V_L \cup \{v_1\}, V_R\cup\{v_2, x\}) = Rs(G, V_L, V_R\cup \{x\})$. Using the recursive definition of the score,  $Rs(P)\leq Ls(G, V_L, V_R\cup \{x\})$, which leads to $Rs(P)\leq Rs(G,V_L\cup\{v_1\},V_R\cup\{v_2\})$.

Finally, we have proved $Rs(P) = Rs(G, V_L \cup \{v_1\}, V_R \cup \{v_2\})$.
\end{proof}

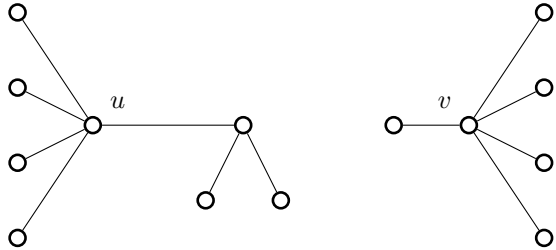
\begin{figure}[ht]
    \centering

\begin{tikzpicture}

\draw (0,0) node[v] (A)  {} node[above right = .15cm] {$u$};
\draw (2,0) node[v] (B) {};
\draw (5,0) node[v] (C) {} node[above left= .15cm]{$v$};
\draw (4,0) node[v] (D) {};

\node[v] (A1) at (-1,-1.5) {};
\node[v] (A2) at (-1,-.5) {};
\node[v] (A3) at (-1,.5) {};
\node[v] (A4) at (-1,1.5) {};

\path[draw, -] (A) to (A1);
\path[draw, -] (A) to (A2);
\path[draw, -] (A) to (A3);
\path[draw, -] (A) to (A4);
\path[draw, -] (A) to (B);

\node[v] (B1) at (1.5,-1) {};
\node[v] (B2) at (2.5,-1) {};

\path[draw, -] (B) to (B1);
\path[draw, -] (B) to (B2);

\node[v] (C1) at (6,-1.5) {};
\node[v] (C2) at (6,-.5) {};
\node[v] (C3) at (6,.5) {};
\node[v] (C4) at (6,1.5) {};

\path[draw, -] (C) to (C1);
\path[draw, -] (C) to (C2);
\path[draw, -] (C) to (C3);
\path[draw, -] (C) to (C4);
\path[draw, -] (C) to (D);

\end{tikzpicture}

    \caption{A graph $G$ for which $Ls(G)=4$ with unique optimal move $u$ and $Rs(G)=2$ with unique optimal move $v$.}
    \label{fig:move_LR_dif}
\end{figure}

Note that, this result is only true for equivalent vertices. In general, a good move for Left is not necessarily a good move for Right. For instance, in Figure~\ref{fig:move_LR_dif}, if Left starts by playing $u$, the score is $4$, and if she starts by playing any other vertex, the score is at most $3$, thus her only optimal move is $u$. If Right starts by playing $v$, the score is $2$, but if he starts by playing any other vertex, the score is at least $3$. Hence, his only optimal move is $v$.

Lemma \ref{twin vertices} is actually very useful to deal with similar vertices. We illustrate its power by computing the score for complete binary trees. A complete binary tree of depth $k$ is a rooted tree such that each vertex at depth $j<k$ has exactly two children (and by definition of the depth, each vertex at depth $k$ is a leaf).

\begin{corollary}
Let $T_k$ be a complete binary tree of depth $k\ge 1$. The scores in Maker-Breaker {\inc} are $Ls(T_k) = 2^{k-1}$ and $Rs(T_k) = 2^{k-1}-1$. 
\end{corollary}

\begin{proof}
Let $T_k$ be a complete binary tree of depth $k$. Its leaves are pairwise equivalent. By Lemma \ref{twin vertices}, we can assume that one leaf other two is given to each player. Then, their parents are pairwise equivalent since the unique free vertex there are adjacent is their father and they are all adjacent to exactly one vertex in $V_L$. 
Thus we can again apply Lemma \ref{twin vertices} and attribute one vertex of depth $k-1$ other two to each player. Going on this reasoning until we each the root, for any pair of vertices having the same parent, Maker and Breaker both get one of them. Then the root is claimed by the first player. Finally, the number of edges taken by Maker satisfies $Ls(T_k) = Ls (T_{k-1}) + Rs(T_{k-1}) + 1$ and $Rs(T_k) = Ls (T_{k-1}) + Rs(T_{k-1})$. Since $Ls(T_0) = Rs(T_0) = 0$, we obtain by induction the result. 
\end{proof}

\section{Complexity of {\sc Maker-Breaker Incidence}} \label{pspace}

In this section, we first prove that \mbinc{} is {\sf PSPACE}-complete. Then, we consider the parameterized complexity of \mbinc{} and prove that it is fixed parameter tractable when parameterized by the neighborhood diversity.

\subsection{{\sc Maker-Breaker Incidence} is {\sf PSPACE}-complete}

Reductions in (positional games) are often made from {\sc POS CNF} (see for example \cite{schaefer1978, reisch1981, stockmeyer1973, stockmeyer1976}). In our cases, we need to deal with scores and not only a structure. To handle this problem, we use a quantified version of {\sc Max-2-SAT} that we proved to be {\sf PSPACE}-complete using {\sc 3-QBF}.

\medskip

\noindent \decisionpb{{\sc Q-Max-2-SAT}}{A quantified boolean formula on the form $\varphi= Q_1x_1, \dots, Q_n x_n, \psi(x_1,x_2,...x_n)$, with $Q_i\in \{\forall,\exists\}$ and $\psi$ a $2$-CNF formula on $x_1,...,x_n$, an integer $k$}{True if at least $k$ clauses of the formula are satisfied. False otherwise.}{1}\\[0.5em]

\noindent \decisionpb{{\sc 3-QBF}}{A quantified boolean formula $\Phi= Q_1x_1, \dots, Q_n x_n, \psi(x_1,x_2,...x_n)$, with $Q_i\in \{\forall,\exists\}$ and $\psi$ a $3$-CNF formula on $x_1,...,x_n$}{True iff $\Phi$ is true.}{1}\\[0.5em]

\begin{theorem}
{\sc Q-Max-2-SAT} is {\sf PSPACE}-complete.
\end{theorem}

\begin{proof}
The proof of {\sf PSPACE}-completeness of {\sc Q-Max-2-SAT} is similar to the proof of {\sf NP}-completeness of {\sc Max-2-SAT} from Papadimitriou~\cite{papadimitriou1994}.

First, {\sc Q-Max-2-SAT} is in {\sf PSPACE}, as any valuation can be computed in polynomial space. Therefore, by a min-max argument, it is possible to compute the number of satisfied clauses in polynomial space.

We provide a reduction from 3-QBF. Let $\phi = Q_1x_1, \dots, Q_n x_n \quad \psi(x_1,x_2,...x_n)$ be a 3-QBF formula on $m$ clauses. For each clause $c_i = l_1^i \vee l_2^i \vee l_3^i$ of $\psi$, we introduce a new variable $d_i$ and construct a set $\mathcal C_i$ of 10 clauses $C_i^1, \dots, C_i^{10}$ of at most 2: 
$$\mathcal C_i= \{(l_1), (l_2), (l_3), (d_i), (\neg l_1 \vee  \neg l_2), (\neg l_1 \vee \neg l_3), (\neg l_2 \vee \neg l_3), (\neg d_1 \vee l_1), (\neg d_1 \vee l_2), (\neg d_1 \vee l_3)\}$$

\begin{claim}
Given any valuation of the literals $l_i$'s, if $c_i$ is satisfied, then there exists a valuation of $d_i$ such that exactly seven clauses in $\mathcal C_i$ are satisfied. Otherwise, at most six clauses of $\mathcal C_i$ are satisfied for any valuation of $d_i$
\end{claim}
  
  \begin{claimproof}
 The proof of the claim is a case analysis depending on the number of literals $l_i$ that are true in $c_i$ (since the literals play a symmetric role). The following tabular gives the number $N_C$ of clauses in $\mathcal C_i$ that are satisfied depending on the number $N_L$ of literals $l_i$ that are true and the valuation of $d_i$. 
 
 \begin{center}
 \begin{tabular}{|c||c|c|c|c|c|c|c|c|}
\hline $N_L$ & 0 & 0 & 1 & 1 & 2 & 2 & 3 & 3 \\ \hline
$d_i$ & F & T & F & T & F & T & F & T \\ \hline
$N_C$ & 6 & 4 & 7 & 6 & 7 & 7 & 6 & 7 \\ \hline
\end{tabular}  
\end{center}

\end{claimproof}

Let $\varphi = Q_1 x_1, \dots, Q_n x_n,\exists d_1, \dots, \exists d_n, \underset{i=1}{\overset{m}{\bigwedge}} \underset{j=1}{\overset{10}{\bigwedge}} C_i^j$ and let $k=7m$.

If $\phi$ is true, then, for any valuation obtained by the $Q_i$'s that makes $\psi$ true, there exists a valuation for each $d_j$ such that there are exactly seven clauses satisfied in each set $\mathcal C_j$. Thus, by taking this valuation for each $d_j$, we have that $k = 7m$ clauses satisfied in $\varphi$.

Reciprocally, if $\phi$ is false, then for any valuation provided by the $Q_i$s, there exists a clause $C_j$ that is not satisfied. Therefore, at most six clauses in $\mathcal C_j$ are satisfied. For the other clauses, at most seven of them are satisfied. Thus the total number of satisfied clauses in $\varphi$ is at most $7m-1 = k-1$.

Finally, the formula $\varphi$ of {\sc Q-Max-2-SAT} has at least $7m$ clauses satisfied if and only if $\phi$ is True. 

Up to add a variable in all the clauses of size 1 and quantifying it with a $\forall$, we can suppose that all the clauses of $\varphi$ have size $2$.
\end{proof}

We now turn to the main proof of this section - that is the proof of the complexity of \mbinc.
\medskip

\noindent \decisionpb{\mbinc}{A graph $G$, an integer $k$, a player $P\in \{\text{Left},\text{Right}\}$.}{True iff the $P$ score of $G$ is at least $k$.}{1}\\[0.5em]

\begin{theorem}\label{thm:pspaceincidence}
\mbinc{} is {\sf PSPACE}-complete.
\end{theorem}

The construction provided in the proof will require some tools to order the moves of both player. Let $P=(G,V_L,V_R)$ be a game position of \inc. Let $u$ and $v$ be free vertices. We say that $v$ {\em dominates} $u$ in $P$ and write $v\geq_P u$ if in any position obtained from $P$, it is always more interesting to play $v$ than $u$. 
More formally, $v\geq_P u$ if for any $V'_L, V'_R$ such that $V_L \subset V'_L$ and $V_R \subset V'_R$, $V'_L\cap V'_R=\emptyset$ and $u,v \notin V'_L \cup V'_R$, we have $Rs(G,V'_L\cup \{u\},V'_R) \ge Rs(G,V'_L\cup \{v\},V'_R)$ and $Ls(G,V'_L\cup \{u\},V'_R) \le Ls(G,V'_L\cup \{v\},V'_R)$.

\begin{lemma} \label{lemma:greaterequalvertices}
Let $G = (V,E)$ be a graph and $P=(G,V_L,V_R)$ a position of \mbinc. Let $u,v$ be two free vertices such that $|N(v) \cap V_L| \ge |N(u) \cap V_L| + |N(u) \setminus N(v) \cap V_F|$. Then $v\geq_P u$.
\end{lemma}

\begin{proof}
Let $\mathcal{S}$ be a strategy in $(G, V_L, V_R)$ that plays $u$ before $v$. We define a strategy $\mathcal{S}'$ that plays $v$ before $u$ as follows:

\begin{itemize}
    \item While $\mathcal{S}$ wants to claim a vertex $w \neq u$, claim $w$.
    \item If $\mathcal{S}$ wants to claim $u$ while $v$ is unclaimed, claim $v$ instead, and still consider that $u$ is claimed in $\mathcal{S}$.
    \item When $\mathcal{S}$ wants to claim $v$, if it is already claimed, claim $u$ instead. If the opponent has claimed $u$, consider that he has claimed $v$, and continue to follow $\mathcal{S}$.
\end{itemize}

Following this strategy, according to the moves of the opponent, all the vertices claimed by $\mathcal{S}$ are claimed by $\mathcal{S}'$, with only a difference on $u$ and $v$ if they are not claimed by the same player.

If $\mathcal{S}$ was a strategy for Left, by following $\mathcal{S'}$, each edge that does not contain $u$ nor $v$ that was claimed by $\mathcal{S}$ is claimed by $\mathcal{S'}$, and reciprocally. Concerning the edges containing $u$ or $v$, Left has scored at most $|N(u) \cap V_L| + |N(u) \cap V_F|$ points on them with $\mathcal{S}$ and $|N(v) \cap V_L| + |N(v) \cap V_F|$ by following $\mathcal{S}'$. Therefore, as $|N(v) \cap V_L| \ge |N(u) \cap V_L| + |N(u) \setminus N(v) \cap V_F|$, Left has score at least the same number of edges following $\mathcal{S}'$. 

The same argument shows that Right will have more edges with a vertex claimed by him by playing $v$ instead of $u$.
\end{proof}

\begin{proof}[Proof of Theorem \ref{thm:pspaceincidence}]

First, \mbinc{} is in {\sf PSPACE} as the game last at most $|V|$ moves and the score is at most $|E|$. Thus, it can be computed in polynomial space, according to Section~6.1 in \cite{Hearn2009}.

We prove that \mbinc{} is {\sf PSPACE}-complete by a reduction from {\sc Q-Max-2-SAT}.
In this proof, we consider a quantified formula as a two-player game. We first assume that the formula has the form $\exists x_{2n} \forall x_{2n-1} \exists x_{2n-2}, \dots \forall x_1 \hspace{.2cm} \psi$ 
, i.e. that the quantifiers $\exists, \forall$ are alternating and starting with a quantifier $\exists$. This can be done for any quantified formulas by adding some vertices with the desired quantifier that are put in no clause, and thus that does not change the number of clauses that are satisfied. The first player, Satisfier, tries to satisfy the formula by choosing the values of the even variables $x_{2k}$ (i.e. that are quantified by an $\exists$-quantifier) while the second player, Falsifier, tries to spoil the formula and turn it to False by choosing the values of the odd variables $x_{2k-1}$ (i.e. that are quantified by a $\forall$-quantifier). This classical technique to transform a quantified formula into a game has been used for instance by Rahman and Watson \cite{rahman2021} to show the {\sf PSPACE}-completeness of Maker-Breaker positional games.

Denote $\psi = \underset{j=1}{ \overset{m}{\bigwedge}} (l^j_1 \vee l^j_2)$ for $l^j_1$, $l_2^j$ some literals. We build a graph $G = (V,E)$ as follows (see Figure \ref{fig:reduction pspace}):

\begin{itemize}
    \item For each variable $x_i$, we create $6mi+3$ vertices. These vertices induce three stars of center $v_i$, $\overline{v_i}$ and $\widetilde{v_i}$, and with $2mi$ leaves each. We will denote by $V_i$ the set $\{v_i, \overline{v_i}, \widetilde{v_i}\}$.
    \item  We consider a function $f$ defined by $f(x_i) = v_i$ and $f(\neg x_i) = \overline{v_i}$. For each clause $C_j = l^j_1 \vee l^j_2$, we add an edge $e_j = (f(l^j_1), f(l^j_2))$.
\end{itemize}

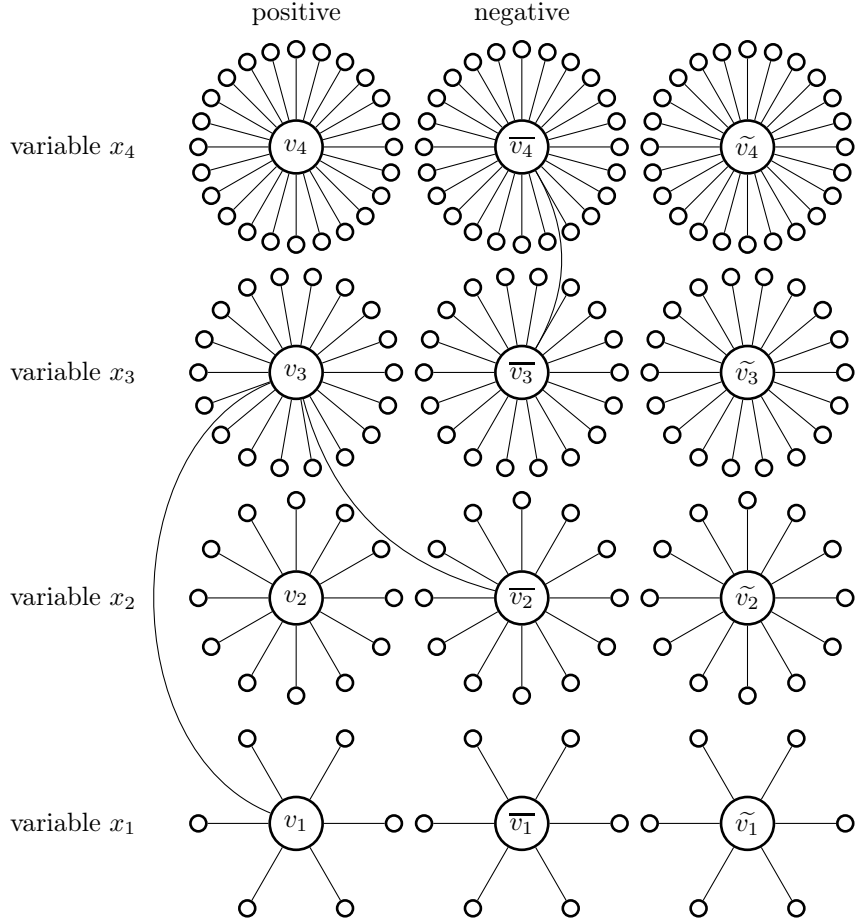
\begin{figure}[ht]
    \centering
\begin{tikzpicture}
\draw ( 1.3 , 0.0 ) node[v](s0){} ;
\draw ( 1.26 , 0.34 ) node[v](s1){} ;
\draw ( 1.13 , 0.65 ) node[v](s2){} ;
\draw ( 0.92 , 0.92 ) node[v](s3){} ;
\draw ( 0.65 , 1.13 ) node[v](s4){} ;
\draw ( 0.34 , 1.26 ) node[v](s5){} ;
\draw ( 0.0 , 1.3 ) node[v](s6){} ;
\draw ( -0.34 , 1.26 ) node[v](s7){} ;
\draw ( -0.65 , 1.13 ) node[v](s8){} ;
\draw ( -0.92 , 0.92 ) node[v](s9){} ;
\draw ( -1.13 , 0.65 ) node[v](s10){} ;
\draw ( -1.26 , 0.34 ) node[v](s11){} ;
\draw ( -1.3 , 0.0 ) node[v](s12){} ;
\draw ( -1.26 , -0.34 ) node[v](s13){} ;
\draw ( -1.13 , -0.65 ) node[v](s14){} ;
\draw ( -0.92 , -0.92 ) node[v](s15){} ;
\draw ( -0.65 , -1.13 ) node[v](s16){} ;
\draw ( -0.34 , -1.26 ) node[v](s17){} ;
\draw ( 0.0 , -1.3 ) node[v](s18){} ;
\draw ( 0.34 , -1.26 ) node[v](s19){} ;
\draw ( 0.65 , -1.13 ) node[v](s20){} ;
\draw ( 0.92 , -0.92 ) node[v](s21){} ;
\draw ( 1.13 , -0.65 ) node[v](s22){} ;
\draw ( 1.26 , -0.34 ) node[v](s23){} ;
\draw ( 0, 0 ) node[c](v4){$v_4$} node[left = 2cm] {variable $x_{4}$} node[above = 1.5cm] {positive};
\path[draw, -] (s0) to (v4);
\path[draw, -] (s1) to (v4);
\path[draw, -] (s2) to (v4);
\path[draw, -] (s3) to (v4);
\path[draw, -] (s4) to (v4);
\path[draw, -] (s5) to (v4);
\path[draw, -] (s6) to (v4);
\path[draw, -] (s7) to (v4);
\path[draw, -] (s8) to (v4);
\path[draw, -] (s9) to (v4);
\path[draw, -] (s10) to (v4);
\path[draw, -] (s11) to (v4);
\path[draw, -] (s12) to (v4);
\path[draw, -] (s13) to (v4);
\path[draw, -] (s14) to (v4);
\path[draw, -] (s15) to (v4);
\path[draw, -] (s16) to (v4);
\path[draw, -] (s17) to (v4);
\path[draw, -] (s18) to (v4);
\path[draw, -] (s19) to (v4);
\path[draw, -] (s20) to (v4);
\path[draw, -] (s21) to (v4);
\path[draw, -] (s22) to (v4);
\path[draw, -] (s23) to (v4);

\draw ( 4.3 , 0.0 ) node[v](s0){} ;
\draw ( 4.26 , 0.34 ) node[v](s1){} ;
\draw ( 4.13 , 0.65 ) node[v](s2){} ;
\draw ( 3.92 , 0.92 ) node[v](s3){} ;
\draw ( 3.65 , 1.13 ) node[v](s4){} ;
\draw ( 3.34 , 1.26 ) node[v](s5){} ;
\draw ( 3.0 , 1.3 ) node[v](s6){} ;
\draw ( 2.66 , 1.26 ) node[v](s7){} ;
\draw ( 2.35 , 1.13 ) node[v](s8){} ;
\draw ( 2.08 , 0.92 ) node[v](s9){} ;
\draw ( 1.87 , 0.65 ) node[v](s10){} ;
\draw ( 1.74 , 0.34 ) node[v](s11){} ;
\draw ( 1.7 , 0.0 ) node[v](s12){} ;
\draw ( 1.74 , -0.34 ) node[v](s13){} ;
\draw ( 1.87 , -0.65 ) node[v](s14){} ;
\draw ( 2.08 , -0.92 ) node[v](s15){} ;
\draw ( 2.35 , -1.13 ) node[v](s16){} ;
\draw ( 2.66 , -1.26 ) node[v](s17){} ;
\draw ( 3.0 , -1.3 ) node[v](s18){} ;
\draw ( 3.34 , -1.26 ) node[v](s19){} ;
\draw ( 3.65 , -1.13 ) node[v](s20){} ;
\draw ( 3.92 , -0.92 ) node[v](s21){} ;
\draw ( 4.13 , -0.65 ) node[v](s22){} ;
\draw ( 4.26 , -0.34 ) node[v](s23){} ;
\draw ( 3, 0 ) node[c](v4b){$\overline{v_4}$} node[above = 1.5cm] {negative};
\path[draw, -] (s0) to (v4b);
\path[draw, -] (s1) to (v4b);
\path[draw, -] (s2) to (v4b);
\path[draw, -] (s3) to (v4b);
\path[draw, -] (s4) to (v4b);
\path[draw, -] (s5) to (v4b);
\path[draw, -] (s6) to (v4b);
\path[draw, -] (s7) to (v4b);
\path[draw, -] (s8) to (v4b);
\path[draw, -] (s9) to (v4b);
\path[draw, -] (s10) to (v4b);
\path[draw, -] (s11) to (v4b);
\path[draw, -] (s12) to (v4b);
\path[draw, -] (s13) to (v4b);
\path[draw, -] (s14) to (v4b);
\path[draw, -] (s15) to (v4b);
\path[draw, -] (s16) to (v4b);
\path[draw, -] (s17) to (v4b);
\path[draw, -] (s18) to (v4b);
\path[draw, -] (s19) to (v4b);
\path[draw, -] (s20) to (v4b);
\path[draw, -] (s21) to (v4b);
\path[draw, -] (s22) to (v4b);
\path[draw, -] (s23) to (v4b);

\draw ( 7.3 , 0.0 ) node[v](s0){} ;
\draw ( 7.26 , 0.34 ) node[v](s1){} ;
\draw ( 7.13 , 0.65 ) node[v](s2){} ;
\draw ( 6.92 , 0.92 ) node[v](s3){} ;
\draw ( 6.65 , 1.13 ) node[v](s4){} ;
\draw ( 6.34 , 1.26 ) node[v](s5){} ;
\draw ( 6.0 , 1.3 ) node[v](s6){} ;
\draw ( 5.66 , 1.26 ) node[v](s7){} ;
\draw ( 5.35 , 1.13 ) node[v](s8){} ;
\draw ( 5.08 , 0.92 ) node[v](s9){} ;
\draw ( 4.87 , 0.65 ) node[v](s10){} ;
\draw ( 4.74 , 0.34 ) node[v](s11){} ;
\draw ( 4.7 , 0.0 ) node[v](s12){} ;
\draw ( 4.74 , -0.34 ) node[v](s13){} ;
\draw ( 4.87 , -0.65 ) node[v](s14){} ;
\draw ( 5.08 , -0.92 ) node[v](s15){} ;
\draw ( 5.35 , -1.13 ) node[v](s16){} ;
\draw ( 5.66 , -1.26 ) node[v](s17){} ;
\draw ( 6.0 , -1.3 ) node[v](s18){} ;
\draw ( 6.34 , -1.26 ) node[v](s19){} ;
\draw ( 6.65 , -1.13 ) node[v](s20){} ;
\draw ( 6.92 , -0.92 ) node[v](s21){} ;
\draw ( 7.13 , -0.65 ) node[v](s22){} ;
\draw ( 7.26 , -0.34 ) node[v](s23){} ;
\draw ( 6, 0 ) node[c](v4t){$\widetilde{v_4}$} ;
\path[draw, -] (s0) to (v4t);
\path[draw, -] (s1) to (v4t);
\path[draw, -] (s2) to (v4t);
\path[draw, -] (s3) to (v4t);
\path[draw, -] (s4) to (v4t);
\path[draw, -] (s5) to (v4t);
\path[draw, -] (s6) to (v4t);
\path[draw, -] (s7) to (v4t);
\path[draw, -] (s8) to (v4t);
\path[draw, -] (s9) to (v4t);
\path[draw, -] (s10) to (v4t);
\path[draw, -] (s11) to (v4t);
\path[draw, -] (s12) to (v4t);
\path[draw, -] (s13) to (v4t);
\path[draw, -] (s14) to (v4t);
\path[draw, -] (s15) to (v4t);
\path[draw, -] (s16) to (v4t);
\path[draw, -] (s17) to (v4t);
\path[draw, -] (s18) to (v4t);
\path[draw, -] (s19) to (v4t);
\path[draw, -] (s20) to (v4t);
\path[draw, -] (s21) to (v4t);
\path[draw, -] (s22) to (v4t);
\path[draw, -] (s23) to (v4t);

\draw ( 1.3 , -3.0 ) node[v](s0){} ;
\draw ( 1.22 , -2.56 ) node[v](s1){} ;
\draw ( 1.0 , -2.17 ) node[v](s2){} ;
\draw ( 0.65 , -1.87 ) node[v](s3){} ;
\draw ( 0.22 , -1.73 ) node[v](s4){} ;
\draw ( -0.22 , -1.73 ) node[v](s5){} ;
\draw ( -0.65 , -1.87 ) node[v](s6){} ;
\draw ( -1.0 , -2.17 ) node[v](s7){} ;
\draw ( -1.22 , -2.56 ) node[v](s8){} ;
\draw ( -1.3 , -3.0 ) node[v](s9){} ;
\draw ( -1.22 , -3.44 ) node[v](s10){} ;
\draw ( -1.0 , -3.83 ) node[v](s11){} ;
\draw ( -0.65 , -4.13 ) node[v](s12){} ;
\draw ( -0.22 , -4.27 ) node[v](s13){} ;
\draw ( 0.22 , -4.27 ) node[v](s14){} ;
\draw ( 0.65 , -4.13 ) node[v](s15){} ;
\draw ( 1.0 , -3.83 ) node[v](s16){} ;
\draw ( 1.22 , -3.44 ) node[v](s17){} ;
\draw ( 0, -3 ) node[c](v3){$v_3$} node[left = 2cm] {variable $x_{3}$};
\path[draw, -] (s0) to (v3);
\path[draw, -] (s1) to (v3);
\path[draw, -] (s2) to (v3);
\path[draw, -] (s3) to (v3);
\path[draw, -] (s4) to (v3);
\path[draw, -] (s5) to (v3);
\path[draw, -] (s6) to (v3);
\path[draw, -] (s7) to (v3);
\path[draw, -] (s8) to (v3);
\path[draw, -] (s9) to (v3);
\path[draw, -] (s10) to (v3);
\path[draw, -] (s11) to (v3);
\path[draw, -] (s12) to (v3);
\path[draw, -] (s13) to (v3);
\path[draw, -] (s14) to (v3);
\path[draw, -] (s15) to (v3);
\path[draw, -] (s16) to (v3);
\path[draw, -] (s17) to (v3);

\draw ( 4.3 , -3.0 ) node[v](s0){} ;
\draw ( 4.22 , -2.56 ) node[v](s1){} ;
\draw ( 4.0 , -2.17 ) node[v](s2){} ;
\draw ( 3.65 , -1.87 ) node[v](s3){} ;
\draw ( 3.22 , -1.73 ) node[v](s4){} ;
\draw ( 2.78 , -1.73 ) node[v](s5){} ;
\draw ( 2.35 , -1.87 ) node[v](s6){} ;
\draw ( 2.0 , -2.17 ) node[v](s7){} ;
\draw ( 1.78 , -2.56 ) node[v](s8){} ;
\draw ( 1.7 , -3.0 ) node[v](s9){} ;
\draw ( 1.78 , -3.44 ) node[v](s10){} ;
\draw ( 2.0 , -3.83 ) node[v](s11){} ;
\draw ( 2.35 , -4.13 ) node[v](s12){} ;
\draw ( 2.78 , -4.27 ) node[v](s13){} ;
\draw ( 3.22 , -4.27 ) node[v](s14){} ;
\draw ( 3.65 , -4.13 ) node[v](s15){} ;
\draw ( 4.0 , -3.83 ) node[v](s16){} ;
\draw ( 4.22 , -3.44 ) node[v](s17){} ;
\draw ( 3, -3 ) node[c](v3b){$\overline{v_3}$} ;
\path[draw, -] (s0) to (v3b);
\path[draw, -] (s1) to (v3b);
\path[draw, -] (s2) to (v3b);
\path[draw, -] (s3) to (v3b);
\path[draw, -] (s4) to (v3b);
\path[draw, -] (s5) to (v3b);
\path[draw, -] (s6) to (v3b);
\path[draw, -] (s7) to (v3b);
\path[draw, -] (s8) to (v3b);
\path[draw, -] (s9) to (v3b);
\path[draw, -] (s10) to (v3b);
\path[draw, -] (s11) to (v3b);
\path[draw, -] (s12) to (v3b);
\path[draw, -] (s13) to (v3b);
\path[draw, -] (s14) to (v3b);
\path[draw, -] (s15) to (v3b);
\path[draw, -] (s16) to (v3b);
\path[draw, -] (s17) to (v3b);

\draw ( 7.3 , -3.0 ) node[v](s0){} ;
\draw ( 7.22 , -2.56 ) node[v](s1){} ;
\draw ( 7.0 , -2.17 ) node[v](s2){} ;
\draw ( 6.65 , -1.87 ) node[v](s3){} ;
\draw ( 6.22 , -1.73 ) node[v](s4){} ;
\draw ( 5.78 , -1.73 ) node[v](s5){} ;
\draw ( 5.35 , -1.87 ) node[v](s6){} ;
\draw ( 5.0 , -2.17 ) node[v](s7){} ;
\draw ( 4.78 , -2.56 ) node[v](s8){} ;
\draw ( 4.7 , -3.0 ) node[v](s9){} ;
\draw ( 4.78 , -3.44 ) node[v](s10){} ;
\draw ( 5.0 , -3.83 ) node[v](s11){} ;
\draw ( 5.35 , -4.13 ) node[v](s12){} ;
\draw ( 5.78 , -4.27 ) node[v](s13){} ;
\draw ( 6.22 , -4.27 ) node[v](s14){} ;
\draw ( 6.65 , -4.13 ) node[v](s15){} ;
\draw ( 7.0 , -3.83 ) node[v](s16){} ;
\draw ( 7.22 , -3.44 ) node[v](s17){} ;
\draw ( 6, -3 ) node[c](v3t){$\widetilde{v_3}$} ;
\path[draw, -] (s0) to (v3t);
\path[draw, -] (s1) to (v3t);
\path[draw, -] (s2) to (v3t);
\path[draw, -] (s3) to (v3t);
\path[draw, -] (s4) to (v3t);
\path[draw, -] (s5) to (v3t);
\path[draw, -] (s6) to (v3t);
\path[draw, -] (s7) to (v3t);
\path[draw, -] (s8) to (v3t);
\path[draw, -] (s9) to (v3t);
\path[draw, -] (s10) to (v3t);
\path[draw, -] (s11) to (v3t);
\path[draw, -] (s12) to (v3t);
\path[draw, -] (s13) to (v3t);
\path[draw, -] (s14) to (v3t);
\path[draw, -] (s15) to (v3t);
\path[draw, -] (s16) to (v3t);
\path[draw, -] (s17) to (v3t);

\draw ( 1.3 , -6.0 ) node[v](s0){} ;
\draw ( 1.13 , -5.35 ) node[v](s1){} ;
\draw ( 0.65 , -4.87 ) node[v](s2){} ;
\draw ( 0.0 , -4.7 ) node[v](s3){} ;
\draw ( -0.65 , -4.87 ) node[v](s4){} ;
\draw ( -1.13 , -5.35 ) node[v](s5){} ;
\draw ( -1.3 , -6.0 ) node[v](s6){} ;
\draw ( -1.13 , -6.65 ) node[v](s7){} ;
\draw ( -0.65 , -7.13 ) node[v](s8){} ;
\draw ( 0.0 , -7.3 ) node[v](s9){} ;
\draw ( 0.65 , -7.13 ) node[v](s10){} ;
\draw ( 1.13 , -6.65 ) node[v](s11){} ;
\draw ( 0, -6 ) node[c](v2){$v_2$} node[left = 2cm] {variable $x_{2}$};
\path[draw, -] (s0) to (v2);
\path[draw, -] (s1) to (v2);
\path[draw, -] (s2) to (v2);
\path[draw, -] (s3) to (v2);
\path[draw, -] (s4) to (v2);
\path[draw, -] (s5) to (v2);
\path[draw, -] (s6) to (v2);
\path[draw, -] (s7) to (v2);
\path[draw, -] (s8) to (v2);
\path[draw, -] (s9) to (v2);
\path[draw, -] (s10) to (v2);
\path[draw, -] (s11) to (v2);
\draw ( 4.3 , -6.0 ) node[v](s0){} ;
\draw ( 4.13 , -5.35 ) node[v](s1){} ;
\draw ( 3.65 , -4.87 ) node[v](s2){} ;
\draw ( 3.0 , -4.7 ) node[v](s3){} ;
\draw ( 2.35 , -4.87 ) node[v](s4){} ;
\draw ( 1.87 , -5.35 ) node[v](s5){} ;
\draw ( 1.7 , -6.0 ) node[v](s6){} ;
\draw ( 1.87 , -6.65 ) node[v](s7){} ;
\draw ( 2.35 , -7.13 ) node[v](s8){} ;
\draw ( 3.0 , -7.3 ) node[v](s9){} ;
\draw ( 3.65 , -7.13 ) node[v](s10){} ;
\draw ( 4.13 , -6.65 ) node[v](s11){} ;
\draw ( 3, -6 ) node[c](v2b){$\overline{v_2}$} ;
\path[draw, -] (s0) to (v2b);
\path[draw, -] (s1) to (v2b);
\path[draw, -] (s2) to (v2b);
\path[draw, -] (s3) to (v2b);
\path[draw, -] (s4) to (v2b);
\path[draw, -] (s5) to (v2b);
\path[draw, -] (s6) to (v2b);
\path[draw, -] (s7) to (v2b);
\path[draw, -] (s8) to (v2b);
\path[draw, -] (s9) to (v2b);
\path[draw, -] (s10) to (v2b);
\path[draw, -] (s11) to (v2b);

\draw ( 7.3 , -6.0 ) node[v](s0){} ;
\draw ( 7.13 , -5.35 ) node[v](s1){} ;
\draw ( 6.65 , -4.87 ) node[v](s2){} ;
\draw ( 6.0 , -4.7 ) node[v](s3){} ;
\draw ( 5.35 , -4.87 ) node[v](s4){} ;
\draw ( 4.87 , -5.35 ) node[v](s5){} ;
\draw ( 4.7 , -6.0 ) node[v](s6){} ;
\draw ( 4.87 , -6.65 ) node[v](s7){} ;
\draw ( 5.35 , -7.13 ) node[v](s8){} ;
\draw ( 6.0 , -7.3 ) node[v](s9){} ;
\draw ( 6.65 , -7.13 ) node[v](s10){} ;
\draw ( 7.13 , -6.65 ) node[v](s11){} ;
\draw ( 6, -6 ) node[c](v2t){$\widetilde{v_2}$} ;
\path[draw, -] (s0) to (v2t);
\path[draw, -] (s1) to (v2t);
\path[draw, -] (s2) to (v2t);
\path[draw, -] (s3) to (v2t);
\path[draw, -] (s4) to (v2t);
\path[draw, -] (s5) to (v2t);
\path[draw, -] (s6) to (v2t);
\path[draw, -] (s7) to (v2t);
\path[draw, -] (s8) to (v2t);
\path[draw, -] (s9) to (v2t);
\path[draw, -] (s10) to (v2t);
\path[draw, -] (s11) to (v2t);

\draw ( 1.3 , -9.0 ) node[v](s0){} ;
\draw ( 0.65 , -7.87 ) node[v](s1){} ;
\draw ( -0.65 , -7.87 ) node[v](s2){} ;
\draw ( -1.3 , -9.0 ) node[v](s3){} ;
\draw ( -0.65 , -10.13 ) node[v](s4){} ;
\draw ( 0.65 , -10.13 ) node[v](s5){} ;
\draw ( 0, -9 ) node[c](v1){$v_1$} node[left = 2cm] {variable $x_{1}$};
\path[draw, -] (s0) to (v1);
\path[draw, -] (s1) to (v1);
\path[draw, -] (s2) to (v1);
\path[draw, -] (s3) to (v1);
\path[draw, -] (s4) to (v1);
\path[draw, -] (s5) to (v1);

\draw ( 4.3 , -9.0 ) node[v](s0){} ;
\draw ( 3.65 , -7.87 ) node[v](s1){} ;
\draw ( 2.35 , -7.87 ) node[v](s2){} ;
\draw ( 1.7 , -9.0 ) node[v](s3){} ;
\draw ( 2.35 , -10.13 ) node[v](s4){} ;
\draw ( 3.65 , -10.13 ) node[v](s5){} ;
\draw ( 3, -9 ) node[c](v1b){$\overline{v_1}$} ;
\path[draw, -] (s0) to (v1b);
\path[draw, -] (s1) to (v1b);
\path[draw, -] (s2) to (v1b);
\path[draw, -] (s3) to (v1b);
\path[draw, -] (s4) to (v1b);
\path[draw, -] (s5) to (v1b);

\draw ( 7.3 , -9.0 ) node[v](s0){} ;
\draw ( 6.65 , -7.87 ) node[v](s1){} ;
\draw ( 5.35 , -7.87 ) node[v](s2){} ;
\draw ( 4.7 , -9.0 ) node[v](s3){} ;
\draw ( 5.35 , -10.13 ) node[v](s4){} ;
\draw ( 6.65 , -10.13 ) node[v](s5){} ;
\draw ( 6, -9 ) node[c](v1t){$\widetilde{v_1}$} ;
\path[draw, -] (s0) to (v1t);
\path[draw, -] (s1) to (v1t);
\path[draw, -] (s2) to (v1t);
\path[draw, -] (s3) to (v1t);
\path[draw, -] (s4) to (v1t);
\path[draw, -] (s5) to (v1t);

\path[draw, -, bend left = 30, ] (v2b) to (v3)  ;
\path[draw, -, bend left = 68, ] (v1) to (v3)  ;
\path[draw, -, bend right = 30, ] (v3b) to (v4b)  ;

\end{tikzpicture}
    \caption{Reduction of $\exists x_4 \forall x_3 \exists x_2 \forall x_1 (\neg x_{2} \vee x_{3}) \wedge (x_{1} \vee x_{3}) \wedge (\neg x_{3} \vee \neg x_{4})$}
    \label{fig:reduction pspace}
\end{figure}

The number of vertices outside sets $V_i$ (i.e. the number of leaves) is $N=\sum_{i=1}^{2n} 6mi= 6mn(2n+1)$. Thus the total number of vertices in $G$ is $N+6n$ and the total number of edges is $N+m$, which is polynomial in the size of $\varphi$. An example of reduction is provided in Figure~\ref{fig:reduction pspace} with $m = 3$ and $n=2$.

Consider a game of {\mbinc} on $G$ with Right starting. Using Lemma~\ref{twin vertices},  for every $1 \le i \le 2n$, the leaves connected to vertices $v_i$, $\overline{v_i}$ and $\widetilde{v_i}$ respectively, are equivalent. Thus, half of them can be given to Left and the other half to Right. Since there are an even number of leaves for each star, the only free vertices after this operation are the  $6n$ vertices in sets $V_i$ for $1 \le i \le n$. Let $P^0=(G,V^0_L,V^0_R)$ be this position, and denote by $V^0_F$ the set of free vertices in this position. By Lemma~\ref{twin vertices}, we have $Rs(G) = Rs(P_0)$.

Now, if $1\leq j <i \leq 2n$, for any $v_i^* \in V_i$ and $v_j^* \in V_j$, we have $|N(v_i^*) \cap V^0_L| = mi$, $|N(v_j^*) \cap V^0_L| = mj$ and  $|N(v_j^*)\cap V^0_F| \le m$. Therefore, by Lemma~\ref{lemma:greaterequalvertices} we have $v^*_i \ge_{P_0} v^*_j$. Moreover, as $N(\widetilde{v_i}) \cap V^0_F = \emptyset$, we also have $v_i \ge_{P_0} \widetilde{v_i}$ and $\overline{v_i} \ge_{P_0} \widetilde{v_i}$. 

Hence, in any optimal strategy in $P_0$ with Right starting, the vertices are played in $n$ rounds, from round $\ell = n$ to $\ell =1$, with the following six steps in each round:

\begin{enumerate}
    \item One vertex chosen by Right among $\{v_{2\ell}, \overline{v_{2\ell}}\}$
    \item The other vertex among $\{v_{2\ell}, \overline{v_{2\ell}}\}$ is taken by Left.
    \item The vertex $\widetilde{v_{2\ell}}$ is taken by Right.
    \item One vertex among $\{v_{2\ell-1}, \overline{v_{2\ell-1}}\}$ is taken by Left.
    \item The second vertex in $\{v_{2\ell-1}, \overline{v_{2\ell-1}}\}$ is taken by Right.
    \item The vertex $\widetilde{v_{2\ell-1}}$ is taken by Left.
\end{enumerate}   

This way, Left will obtain exactly $N' = \sum_{\ell=1}^n (2\ell m + 2(2\ell-1)m)=3mn(n+1)-2mn$ edges in the stars and maybe some other edges in the clause edges. 
Let $k' = N' + m - k+1 $. 

We will prove that $Rs(G)\geq k'$ at \mbinc{} if and only if Falsifier wins at {\sc Q-Max-2-SAT}  on $(\varphi,k)$.

\medskip

\begin{claim}
If Satisfier has a strategy to satisfy $k$ clauses in $\varphi$, then $Rs(G) < k'$.
\end{claim}

\begin{claimproof}
We suppose that Satisfier has a winning strategy $\mathcal{S}$ in $(\varphi, k)$. 
We consider that both Right and Left play optimally in $G$ and thus we can assume that the game is played in $P_0$ and respects the previous order. 

Consider the following strategy for Right. At each round $\ell$ from $\ell=n$ to $\ell=1$, Right takes a decision only at Step 1. If Satisfier would turn $x_{2i}$ to True in the game played on $\varphi$, then Right plays $v_{2i}$, otherwise, he plays $\overline{v_{2i}}$. Then, Steps 2 and 3 are determined. At Step 4, if Left plays $v_{2i-1}$ then Right considers that Falsifier has turned $x_{2i-1}$ to False, otherwise he considers she has turned it to True. Then again, Steps 5 and 6 are determined.
By following this strategy, the underlying value obtained for $\varphi$ is exactly the value that Satisfier would obtain by playing according to $\mathcal{S}$. Thus, at least $k$ clauses are satisfied in $\varphi$. 

Note that for a literal $l^j$, the vertex $f(l^j)$ is taken by Right if and only if $l^j$ is True in the game of {\sc Q-Max-2-SAT}. Let $C_j = l^j_1 \vee l^j_2$ be a clause. If Left has claimed the two extremities of $e_j$, it means that Left has played $f(l^j_1)$ and $f(l^j_2)$. Therefore, the underlying values of $l^j_1$ and of $l^j_2$ are both False, and $C_j$ is not satisfied in $\psi$. Hence, Left claims at most $m- k$ edges $e_j$. Finally, Left claimed at most $k'-1$ edges and we have $Rs(G) < k'$. 
\end{claimproof}

\medskip

\begin{claim}
If Falsifier has a strategy such that at most $k-1$ clauses are satisfied in $\phi$, then $Rs(G) \ge k'$.
\end{claim}

\begin{claimproof}
We now suppose that Falsifier has a winning strategy $\mathcal{S}$ in $(\varphi, k)$. 
We consider that both Right and Left play optimally in $G$ and thus we can assume that the game is played in $P_0$ and respects the previous order. Consider the following strategy for Left. At each round $\ell$ from $\ell=n$ to $\ell=1$, Left takes a decision only at Step 4. At Step 1, if Right plays $v_{2\ell}$ then Left considers that Satisfier has turned $x_{2\ell}$ to True, otherwise she considers he has turned it to False. Then, Steps 2 and 3 are determined. At Step 4, if Falsifier would turn $x_{2i-1}$ to False in the game played on $\varphi$, then Left plays $v_{2i-1}$, otherwise, she plays $\overline{v_{2i-1}}$. Then again, Steps 5 and 6 are determined.

By following this strategy, the underlying value obtained for $\varphi$ is exactly the value that Falsifier would obtain by playing according to $\mathcal{S}$. Thus, it would satisfy at most $k-1$ clauses in $\varphi$.
As before, if a clause $l^j_1\vee l^j_2$ is not satisfied in $\varphi$ it means that both vertices $f(l^j_1)$ and $f(l^j_2)$ are taken by Left and thus Left got the edge.
Thus Left claims at least $N'+m-k+1$ edges in the game $G$ and $Rs(G) \ge k'$. 
\end{claimproof}
\end{proof}

\begin{remark}
Note that, up to add a useless variable in $\varphi$, $\varphi$ could start by a $\forall$-quantifier, implying that \mbinc{} is {\sf PSPACE}-complete even if Left starts.
\end{remark}

\begin{corollary}
$3$-uniform {\MMSPG} is {\sf PSPACE}-complete.
\end{corollary}

\begin{proof}
The proof is similar to the second part of the proof of Corollary \ref{cor:generalcomplexity}. From a graph $G = (V,E)$ of \mbinc, we consider the instance of $3$-uniform \MMSPG{} obtained by adding a universal vertex $v_0$. Consider the hypergraph $H =(V \cup \{v_0\}, \left \{e \cup \{v_0\} | e \in E \right \})$. When Left starts, any optimal strategy starts by playing $v_0$, otherwise Right plays it and the score will be at most $0$. Then we are left to a Maker-Breaker position as Right cannot score any point, but starts. Finally the Left score of $H$ in Maker-Maker convention is equal to the Right score of $G$ in Maker-Breaker convention, which is {\sf PSPACE}-complete to compute.
\end{proof}

\subsection{Complexity parameterized by the neighborhood diversity}

Neighborhood diversity is a graph parameter introduced by Lampis \cite{lampis2010} to generalize FPT algorithms parameterized by vertex cover to larger classes of graphs.
Let $G$ be a graph. We say that two vertices $u$ and $v$ have the same {\em type} if $N(v)\setminus \{u\}=N(u)\setminus \{v\}$. The graph $G$ has {\em neighborhood diversity} at most $w$ if there exists a partition of $V$ into at most $w$ sets such that the vertices in each set have all the same type. Note that each set must induce a clique or an independent set. Furthermore, if a graph has bounded vertex cover, then it has bounded neighborhood diversity.

A decision problem has a {\em kernel} for a parameter $w$, if for any parameterized instance $(P,w)$ of the problem, there exists an instance $(P',w')$ and a computable function $f$, such that $P$ reduces to $P'$ in polynomial time in $(|P|,w)$ and such that $|P'|,|w'| \le f(w)$. If $f(w) = O(w^3)$, the kernel is said to be cubic. If $f(w) = O(w \log(w))$, the kernel is said to be quasilinear. Having a kernel implies that the problem is fixed-parameter tractable for this parameter.  

\begin{theorem}
\mbinc{} parameterized by the neighborhood diversity $w$ has a cubic kernel.
\end{theorem}

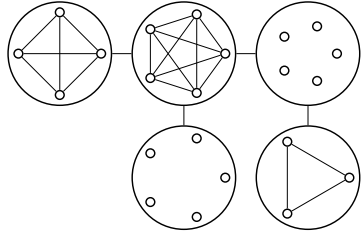
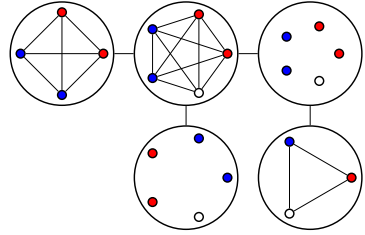
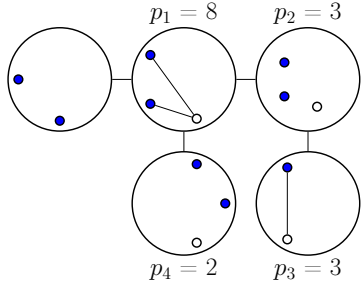
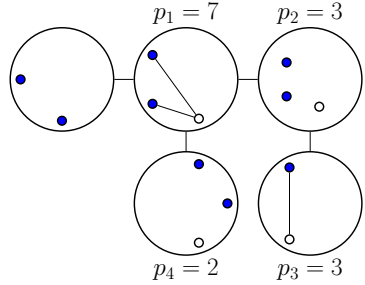
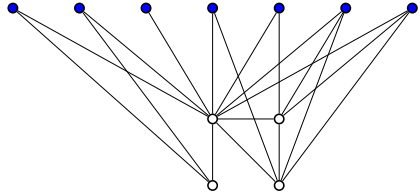
\begin{figure}[ht]
\centering
\begin{subfigure}[b]{.47 \textwidth}
\centering
\scalebox{.55}{
\begin{tikzpicture}

\tikzstyle{gros}=[circle, inner sep=0, minimum size = 2.5cm, line width = 1pt, draw=black, fill=white, text= black]

 \node[gros] (A) at (0,0) {};
 \node[gros] (B) at (3,0) {};
 \node[gros] (C) at (-3,0) {};
 \node[gros] (D) at (0,-3) {};
 \node[gros] (E) at (3,-3) {};

\path[draw, -] (A) to (B);
\path[draw, -] (A) to (C);
\path[draw, -] (A) to (D);
\path[draw, -] (B) to (E);

\draw ( 1.0 , 0.0 ) node[v](s0){} ;
\draw ( 0.31 , 0.95 ) node[v](s1){} ;
\draw ( -0.81 , 0.59 ) node[v](s2){} ;
\draw ( -0.81 , -0.59 ) node[v](s3){} ;
\draw ( 0.31 , -0.95 ) node[v](s4){} ;
\path[draw, -] (s0) to (s1);
\path[draw, -] (s0) to (s2);
\path[draw, -] (s0) to (s3);
\path[draw, -] (s0) to (s4);
\path[draw, -] (s1) to (s2);
\path[draw, -] (s1) to (s3);
\path[draw, -] (s1) to (s4);
\path[draw, -] (s2) to (s3);
\path[draw, -] (s2) to (s4);
\path[draw, -] (s3) to (s4);

\draw ( 3.7 , 0.0 ) node[v](s0){} ;
\draw ( 3.22 , 0.66 ) node[v](s1){} ;
\draw ( 2.43 , 0.41 ) node[v](s2){} ;
\draw ( 2.43 , -0.41 ) node[v](s3){} ;
\draw ( 3.22 , -0.66 ) node[v](s4){} ;

\draw ( -2.0 , 0.0 ) node[v](s0){} ;
\draw ( -3.0 , 1.0 ) node[v](s1){} ;
\draw ( -4.0 , 0.0 ) node[v](s2){} ;
\draw ( -3.0 , -1.0 ) node[v](s3){} ;
\path[draw, -] (s0) to (s1);
\path[draw, -] (s0) to (s2);
\path[draw, -] (s0) to (s3);
\path[draw, -] (s1) to (s2);
\path[draw, -] (s1) to (s3);
\path[draw, -] (s2) to (s3);

\draw ( 1.0 , -3.0 ) node[v](s0){} ;
\draw ( 0.31 , -2.05 ) node[v](s1){} ;
\draw ( -0.81 , -2.41 ) node[v](s2){} ;
\draw ( -0.81 , -3.59 ) node[v](s3){} ;
\draw ( 0.31 , -3.95 ) node[v](s4){} ;

\draw ( 4.0 , -3.0 ) node[v](s0){} ;
\draw ( 2.5 , -2.13 ) node[v](s1){} ;
\draw ( 2.5 , -3.87 ) node[v](s2){} ;
\path[draw, -] (s0) to (s1);
\path[draw, -] (s1) to (s2);
\path[draw, -] (s2) to (s0);

\end{tikzpicture}}
\caption{A graph to kernelize.  Set $k = 30$.}
\label{FPT:STEP0}
\end{subfigure} 
\begin{subfigure}[b]{.47 \textwidth}
\centering
\scalebox{.55}{

\begin{tikzpicture}
\tikzstyle{gros}=[circle, inner sep=0, minimum size = 2.5cm, line width = 1pt, draw=black, fill=white, text= black]

 \node[gros] (A) at (0,0) {};
 \node[gros] (B) at (3,0) {};
 \node[gros] (C) at (-3,0) {};
 \node[gros] (D) at (0,-3) {};
 \node[gros] (E) at (3,-3) {};

\path[draw, -] (A) to (B);
\path[draw, -] (A) to (C);
\path[draw, -] (A) to (D);
\path[draw, -] (B) to (E);

\draw ( 1.0 , 0.0 ) node[R](s0){} ;
\draw ( 0.31 , 0.95 ) node[R](s1){} ;
\draw ( -0.81 , 0.59 ) node[L](s2){} ;
\draw ( -0.81 , -0.59 ) node[L](s3){} ;
\draw ( 0.31 , -0.95 ) node[v](s4){} ;
\path[draw, -] (s0) to (s1);
\path[draw, -] (s0) to (s2);
\path[draw, -] (s0) to (s3);
\path[draw, -] (s0) to (s4);
\path[draw, -] (s1) to (s2);
\path[draw, -] (s1) to (s3);
\path[draw, -] (s1) to (s4);
\path[draw, -] (s2) to (s3);
\path[draw, -] (s2) to (s4);
\path[draw, -] (s3) to (s4);

\draw ( 3.7 , 0.0 ) node[R](s0){} ;
\draw ( 3.22 , 0.66 ) node[R](s1){} ;
\draw ( 2.43 , 0.41 ) node[L](s2){} ;
\draw ( 2.43 , -0.41 ) node[L](s3){} ;
\draw ( 3.22 , -0.66 ) node[v](s4){} ;

\draw ( -2.0 , 0.0 ) node[R](s0){} ;
\draw ( -3.0 , 1.0 ) node[R](s1){} ;
\draw ( -4.0 , 0.0 ) node[L](s2){} ;
\draw ( -3.0 , -1.0 ) node[L](s3){} ;
\path[draw, -] (s0) to (s1);
\path[draw, -] (s0) to (s2);
\path[draw, -] (s0) to (s3);
\path[draw, -] (s1) to (s2);
\path[draw, -] (s1) to (s3);
\path[draw, -] (s2) to (s3);

\draw ( 1.0 , -3.0 ) node[L](s0){} ;
\draw ( 0.31 , -2.05 ) node[L](s1){} ;
\draw ( -0.81 , -2.41 ) node[R](s2){} ;
\draw ( -0.81 , -3.59 ) node[R](s3){} ;
\draw ( 0.31 , -3.95 ) node[v](s4){} ;

\draw ( 4.0 , -3.0 ) node[R](s0){} ;
\draw ( 2.5 , -2.13 ) node[L](s1){} ;
\draw ( 2.5 , -3.87 ) node[v](s2){} ;
\path[draw, -] (s0) to (s1);
\path[draw, -] (s1) to (s2);
\path[draw, -] (s2) to (s0);

\end{tikzpicture}}
\caption{Step~1, $k = 30$.} 
\label{FPT:STEP1}
\end{subfigure}

\vspace{1cm}

\begin{subfigure}[b]{.47 \textwidth}
    \centering
\scalebox{.55}{
\begin{tikzpicture}
\tikzstyle{gros}=[circle, inner sep=0, minimum size = 2.5cm, line width = 1pt, draw=black, fill=white, text= black]

 \node[gros] (A) at (0,0) {};
 \node[gros] (B) at (3,0) {};
 \node[gros] (C) at (-3,0) {};
 \node[gros] (D) at (0,-3) {};
 \node[gros] (E) at (3,-3) {};
 \node[font = \LARGE] at (0,1.6) { $p_1$ = 8};
 \node[font = \LARGE] at (3,1.6) { $p_2$ = 3};
 \node[font = \LARGE] at (0,-4.6) { $p_4$ = 2};
 \node[font = \LARGE] at (3,-4.6) { $p_3$ = 3};

\path[draw, -] (A) to (B);
\path[draw, -] (A) to (C);
\path[draw, -] (A) to (D);
\path[draw, -] (B) to (E);


\draw ( -0.81 , 0.59 ) node[L](s2){} ;
\draw ( -0.81 , -0.59 ) node[L](s3){} ;
\draw ( 0.31 , -0.95 ) node[v](s4){} ;

\path[draw, -] (s2) to (s4);
\path[draw, -] (s3) to (s4);

\draw ( 2.43 , 0.41 ) node[L](s2){} ;
\draw ( 2.43 , -0.41 ) node[L](s3){} ;
\draw ( 3.22 , -0.66 ) node[v](s4){} ;


\draw ( -4.0 , 0.0 ) node[L](s2){} ;
\draw ( -3.0 , -1.0 ) node[L](s3){} ;

\draw ( 1.0 , -3.0 ) node[L](s0){} ;
\draw ( 0.31 , -2.05 ) node[L](s1){} ;
\draw ( 0.31 , -3.95 ) node[v](s4){} ;

\draw ( 2.5 , -2.13 ) node[L](s1){} ;
\draw ( 2.5 , -3.87 ) node[v](s2){} ;
\path[draw, -] (s1) to (s2);

\end{tikzpicture}}
\caption{Step~2, 16 edges removed, $k = 30 - 16 = 14$}
\label{FPT:STEP2}
\end{subfigure}
\begin{subfigure}[b]{.47 \textwidth}
\centering
    \scalebox{.55}{
\begin{tikzpicture}
\tikzstyle{gros}=[circle, inner sep=0, minimum size = 2.5cm, line width = 1pt, draw=black, fill=white, text= black]

 \node[gros] (A) at (0,0) {};
 \node[gros] (B) at (3,0) {};
 \node[gros] (C) at (-3,0) {};
 \node[gros] (D) at (0,-3) {};
 \node[gros] (E) at (3,-3) {};
 \node[font = \LARGE] at (0,1.6) { $p_1 = 7$};
 \node[font = \LARGE] at (3,1.6) { $p_2 = 3$};
 \node[font = \LARGE] at (0,-4.6) { $p_4 = 2$};
 \node[font = \LARGE] at (3,-4.6) { $p_3 = 3$};

\path[draw, -] (A) to (B);
\path[draw, -] (A) to (C);
\path[draw, -] (A) to (D);
\path[draw, -] (B) to (E);


\draw ( -0.81 , 0.59 ) node[L](s2){} ;
\draw ( -0.81 , -0.59 ) node[L](s3){} ;
\draw ( 0.31 , -0.95 ) node[v](s4){} ;

\path[draw, -] (s2) to (s4);
\path[draw, -] (s3) to (s4);

\draw ( 2.43 , 0.41 ) node[L](s2){} ;
\draw ( 2.43 , -0.41 ) node[L](s3){} ;
\draw ( 3.22 , -0.66 ) node[v](s4){} ;


\draw ( -4.0 , 0.0 ) node[L](s2){} ;
\draw ( -3.0 , -1.0 ) node[L](s3){} ;

\draw ( 1.0 , -3.0 ) node[L](s0){} ;
\draw ( 0.31 , -2.05 ) node[L](s1){} ;
\draw ( 0.31 , -3.95 ) node[v](s4){} ;

\draw ( 2.5 , -2.13 ) node[L](s1){} ;
\draw ( 2.5 , -3.87 ) node[v](s2){} ;
\path[draw, -] (s1) to (s2);

\end{tikzpicture}}
\caption{Step~3, $p_1$ has decreased by $1$. $k = 13$.}
\label{FPT:STEP3}
\end{subfigure}

\vspace{1cm}

\begin{subfigure}[b]{.47 \textwidth}
    \begin{center}
    \scalebox{.59}{
\begin{tikzpicture}
\tikzstyle{gros}=[circle, inner sep=0, minimum size = 2.5cm, line width = 1pt, draw=black, fill=white, text= black]

 \node[v] (A) at (0,0) {};
 \node[v] (B) at (1.5,0) {};
 \node[v] (D) at (0,-1.5) {};
 \node[v] (E) at (1.5,-1.5) {};

 \node[L] (L1) at (-4.5,2.5) {};
 \node[L] (L2) at (-3,2.5) {};
 \node[L] (L3) at (-1.5,2.5) {};
 \node[L] (L4) at (0,2.5) {};
 \node[L] (L5) at (1.5,2.5) {};
 \node[L] (L6) at (3,2.5) {};
 \node[L] (L7) at (4.5,2.5) {};

\path[draw, -] (A) to (L1);
\path[draw, -] (A) to (L2);
\path[draw, -] (A) to (L3);
\path[draw, -] (A) to (L4);
\path[draw, -] (A) to (L5);
\path[draw, -] (A) to (L6);
\path[draw, -] (A) to (L7);

\path[draw, -] (B) to (L5);
\path[draw, -] (B) to (L6);
\path[draw, -] (B) to (L7);

\path[draw, -] (E) to (L4);
\path[draw, -] (E) to (L6);
\path[draw, -] (E) to (L7);

\path[draw, -] (D) to (L1);
\path[draw, -] (D) to (L2);

\path[draw, -] (A) to (B);
\path[draw, -] (A) to (D);
\path[draw, -] (A) to (E);
\path[draw, -] (B) to (E);

\end{tikzpicture}}
\end{center}
\caption{Step~4, each vertex $v_i$ has $p_i$ blue neighbors.}
\label{FPT:STEP4}
    \end{subfigure}
    \caption{Example of a kernelization. Vertices in the same circle have same type. An edge between two circles means that all the edges between the vertices of the two circles are in the graph. Blue and red vertices are given to Left and Right respectively. We start with $n = 22$ and after Step~1 $r=4$.}
    \label{fig:FPT}
\end{figure}

\begin{proof}
In this proof, we will consider as instances of \mbinc{} triplets $(P,k,Left)$ where $P$ is a position of \mbinc{} played on $G$ (i.e. some vertices are already played). Note that this does not change the complexity of the problem. Indeed, from any position $P=(G,V_L,V_R)$ one can obtain a graph $G'$ with no vertices played for which the games are equivalent. First remove all the vertices in $V_R$ of the graph. Then, duplicate each vertex in $V_L$ by creating a twin vertex having the same neighbourhood and free the vertices in $V_L$. By Lemma \ref{twin vertices}, one can assume that both players will take one vertex in each pair of twins.

Let $G=(V,E)$ be a graph of neighborhood diversity $w$. 
Consider a partition $(V_1, \dots, V_w)$ of $V$ such that the vertices in each part are all of the same type. We provide the following kernelization algorithm. Let $I=((G,\emptyset,\emptyset),k,P)$ where $P\in \{Left,Right\}$ be an instance of \mbinc{}.  
An example of the different steps is provided in Figure~\ref{fig:FPT}.

{\bf \underline{Step~1}:} While there exists a part $V_i$, $1\le i \le w $ such that there are at least two free vertices $u,v \in V_i$, add $u$ to $V_L$ and $v$ to $V_R$. By Lemma~\ref{twin vertices}, this transformation does not change the outcome of the game. At the end of Step~1, there are at most $w$ free vertices in $G$. In Figure~\ref{fig:FPT}\subref{FPT:STEP1}, it consists in distributing vertices of same type between Left and Right.

{\bf \underline{Step~2}:} Remove all the edges included in $V_L$ and set $k \leftarrow k - |e \subset V_L|$. Then remove from $G$ all the vertices in $V_R$ that cannot count for any point. This transformation do not change the outcome of $I$. At this moment, $G$ only contains free vertices or vertices claimed by Left, and any edge has at least one free extremity. In Figure~\ref{fig:FPT}\subref{FPT:STEP2}, it consists in removing the $16$ edges on which the two endpoints are claimed by Left, and to remove the red vertices and their incident edges. Therefore, $k$ is decreased from $30$ to $14$.

{\bf \underline{Step~3}:} Let $r$ the number of free vertices in $P$, we have $r \le w$. Let $v_1, \dots, v_{r}$ be these vertices. For $1 \le i \le r$, let $p_i = |N(v_i) \cap V_L|$ and order the vertices such that $p_1 \ge p_2 \ge \dots \ge p_{r}$. 
While there exists an integer $i$ such that $p_{i} > p_{i+1} +r$ (with $p_{r+1} = 0$), by Lemma~\ref{lemma:greaterequalvertices}, there exists an optimal strategy in which the vertices $v_1, \dots, v_{i}$ are played before the vertices $v_{i+1}, \dots, v_{r}$. On these vertices, Left will score at least $p_{i}$ at each Left move. Therefore, we can do the following transformation. Let $s=p_{i}-p_{i+1} - r$ for any $1 \le j \le i$, set $p_j \leftarrow p_j - s$  and set $k \leftarrow k - s\left \lceil \frac{i}{2}\right \rceil $. 
Repeat Step~3 until we have $p_i\leq p_{i+1}+r$ for all $1\leq i\leq r$.
In particular, we have after these operations $p_1 \le r^2$. 
In Figure~\ref{fig:FPT}\subref{FPT:STEP3}, it happens only once, as $p_1 = 8$, $p_2 = 3$ and $w' = 4$. Therefore, we set $p_1 = 7$ and $k$ is decreased from $14$ to $13$.

{\bf \underline{Step~4}:} Let $U = \{u_1, \dots, u_{p_1}\}$ be $p_1$ new vertices and transform $(G,V_L,\emptyset)$ into $((G\setminus V_L) \cup U, U, \emptyset)$, and, for $1 \le i \le r$, connect the vertex $v_i$ to any $p_i$ vertices in $U$. This transformation do not change the outcome of the game, since only the number of neighbors in $V_L$ matters when a vertex is played. In Figure~\ref{fig:FPT}\subref{FPT:STEP4}, we have $p_1 = 7$. Thus, $U$ contains seven vertices and each remaining uncolored vertex $v_i$ is connected to $p_i$ of these seven vertices.

Finally, if $k \ge r^3$, as there are at most $r^3$ edges in the final graph, we can just transform $P$ into a trivial False instance like the empty graph with $k = 1$. Thus, we can assume that $k \le r^3$.

The instance obtained has $p_1 + r \le r^2+r\le w^2 + w$  vertices, at most $r * p_1 \le r^3 \le w^3$ edges, $k \le r^3 \le w^3$ and the same outcome as the input. Finally, this new instance has cubic size in $w$ and thus \mbinc{} has a cubic kernel.
\end{proof}

\begin{corollary}
Let $G$ be a graph of order $n$ and neighborhood diversity $w$. In \mbinc{} $Ls(G)$ and $Rs(G)$ can be computed in time $O(w^2w! + n^2)$
\end{corollary}

\begin{proof}
We can compute the kernel in time $n^2$, and then try all the possible games by testing all the moves in time $w^2w!$.
\end{proof}

Note that the cubic size of the kernel is mostly due to the $w^2$ vertices that are already claimed by Left. As these vertices cannot be played any longer, by giving weight to the vertices, it is possible to have a quasilinear kernel by storing only the number of neighbors of each vertex that are already claimed by Left instead of vertices themselves.

\section{Paths and cycles}\label{section path}

We here give the exact values of the score for {\mbinc} played on paths and cycles. For that purpose, we will consider the equivalence properties of Milnor's universe detailed in Section 2. In particular, the notion of negative will be required, implying to consider the partisan version of {\inc}. More precisely, in this section, instances of {\mbinc} will correspond to paths or cycles where the edges are either colored all blue (i.e. only Left can get points) or all red (i.e. only Right can get points). The notations are defined as follows: 
\begin{itemize}

\item $P^L_n$: path of order $n$ where all the edges are colored blue. We denote the vertices of $P^L_n$ by $\{v_0, \dots, v_{n-1}\}$

\item $P^R_n$: path of order $n$ where all the edges are colored red. We denote the vertices of $P^R_n$ by $\{v'_0, \dots, v'_{n-1}\}$

\end{itemize}

By definition, we have that $P^L_n=-P^R_n$. 

\subsection{Equivalences of paths}

We first give the main result about the equivalence between paths modulo $5$. To present it, we introduce a usual notation in scoring game theory: for $k \in \mathbb{Z}$, we define by $k$ the game with no option and where Left has a score of $k$ points. Thus, in \mbinc, the game $1$ is equivalent to $P^L_2$ in which Left has claimed the two vertices and $-1$ is equivalent to $P^R_2$ in which Right has claimed the two vertices. Note that for any game $G$ and any integer $k$, we have $G \equiv k$ if and only if $Ls(G) = Rs(G) = k$. The main theorem of this section states that paths of order at least $6$ are equivalent to paths having five vertices less, with a difference of one in the score. This result remains true if an extremity of the path is already colored.

\begin{theorem}\label{equivalence paths}
Let $n \ge 1$ be an integer. We have $P^L_{n+5} \equiv P^L_n +1$ and $P^R_{n+5} \equiv P^R_n -1$.

\noindent
Let $n \ge 2$ be an integer. We have $(P^L_{n+5}, \{v_0\}, \emptyset) \equiv (P^L_{n}, \{v_0\}, \emptyset) + 1$ and $(P^R_{n+5}, \emptyset, \{v'_0\}) \equiv (P^R_{n}, \emptyset, \{v'_0\}) - 1$.
\end{theorem}

The rest of this subsection will be dedicated to the proof of this theorem.

\subsubsection{Strategy for Left when Right starts}

\begin{lemma}\label{strategy left}
Let $n \ge 1$ be an integer. In \mbinc{}, we have $Rs(P^L_{n+5} + P^R_n) \ge 1$. 

Let $n \ge 2$ be an integer. In \mbinc{}, we have  $Rs(P^L_{n+5} + P^R_n, \{v_0\}, \{v'_0\}) \ge 1$.
\end{lemma}

This proof will be done by induction. Therefore, to handle the small cases, the scores of first paths will be required. They are recorded in Figure~\ref{fig:score_path} and Figure~\ref{fig:score_path_A} and can be easily checked by hand.

\begin{figure}[ht]
    \centering
    \begin{tabular}{|c|c|c|c|c|c|c|c|c|c|c|}
    \hline
    $n$ & $1$ & $2$ &$3$ & $4$ & $5$ & $6$ &$7$& $8$ & $9$ & $10$ \\
    \hline
    $Ls(P^L_n)$ & $0$ & $0$ &$1$ & $1$ & $1$ & $1$ &$1$& $2$ & $2$ & $2$ \\
    \hline
    $Rs(P^L_n)$ & $0$ & $0$ &$0$ & $0$ & $0$ & $1$ &$1$& $1$ & $1$ & $1$ \\
    \hline
\end{tabular}
    \caption{First scores in short paths}
    \label{fig:score_path}
\end{figure}

\begin{figure}[ht]
    \centering
    \begin{tabular}{|c|c|c|c|c|c|c|c|c|c|c|c|}
    \hline
    $n$ & $1$ & $2$ & $3$ &$4$ & $5$ & $6$ & $7$ &$8$& $9$ & $10$ & $11$ \\
    \hline
    $Ls((P^L_n, \{v_0\}, \emptyset))$ & $0$ & $1$ & $1$ &$1$ & $1$ & $2$ & $2$ &$2$& $2$ & $2$ & $3$ \\
    \hline
    $Rs((P^L_n, \{v_0\}, \emptyset))$ & $0$ &$0$ & $0$ &$0$ & $1$ & $1$ & $1$ &$1$& $1$ & $2$ & $2$ \\
    \hline
\end{tabular}
    \caption{First scores in short paths with an extremity claimed by Left}
    \label{fig:score_path_A}
\end{figure}

\begin{proof}
In order to prove that $Rs(P^L_{n+5} + P^R_n) \ge 1$ ($Rs(P^L_{n+5} + P^R_n, \{v_0\}, \{v'_0\}) \ge 1$ resp.), we provide a strategy for Left by induction. If $1 \le n \le 5$ ($2 \le n \le 6$ resp.), a computation can verify that the result is true.

If $n \ge 6$ ($n \ge 7 $ resp.), we consider the first move of Right:
\begin{itemize}
    \item If Right plays a vertex $v'_i$ for $0 \le i \le n-1$ ($1 \le i \le n-1$ resp.), Left answers by playing the vertex $v_i$. The resulting position is $(P^L_{n+5} + P^R_n, \{v_i\}, \{v'_i\})$ ($(P^L_{n+5} + P^R_n, \{v_0, v_i\}, \{v'_0, v'_i\})$  resp.), which is equivalent to $(P^L_{i+1} + P^R_{i+1}, \{v_i\}, \{v'_i\}) + (P^L_{n+5-i} + P^R_{n-i}, \{v_0\}, \{v'_0\})$ ($(P^L_{i+1} + P^R_{i+1}, \{v_0, v_i\}, \{v'_0, v'_i\}) + (P^L_{n+5-i} + P^R_{n-i}, \{v_0\}, \{v'_0\})$ resp.). As we have $(P^L_{i+1} + P^R_{i+1}, \{v_i\}, \{v'_i\}) \equiv 0$  ($(P^L_{i+1} + P^R_{i+1}, \{v_0, v_i\}, \{v'_0, v'_i\}) \equiv 0$ resp.) and $(P^L_{n+5-i} + P^R_{n-i}, \{v_0\}, \{v'_0\})$ satisfies the induction hypothesis, and therefore the score is at least one.
    \item If Right plays a vertex $v_i$ for $0 \le i \le n-1$ ($1 \le i \le n-1$ resp.), Left answers by playing the vertex $v'_i$. The resulting position is $(P^L_{n+5} + P^R_n, \{v'_i\}, \{v_i\})$ ($(P^L_{n+5} + P^R_n, \{v_0, v'_i\}, \{v'_0, v_i\})$  resp.), which is equivalent to $(P^L_{i} + P^R_{i}) + (P^L_{n+5-(i+1)} + P^R_{n-(i+1)})$ ($(P^L_{i} + P^R_{i}, \{v_0\}, \{v'_0\}) + (P^L_{n+5-(i+1)} + P^R_{n-(i+1)})$ resp.). As we have $(P^L_{i} + P^R_{i}) \equiv 0$  ($(P^L_{i} + P^R_{i}, \{v_0\}, \{v'_0\}) \equiv 0$ resp.) and $(P^L_{n+5-(i+1)} + P^R_{n-(i+1)})$ satisfies the induction hypothesis, the score is at least one.
    \item If Right plays a vertex $v_i$ for $n \le i \le n+4$. Left answers by playing $v'_{i-5}$, which exists as $n \ge 6$ ($n \ge 7$ resp.). The resulting position is $(P^L_{n+5} + P^R_n, \{v'_{i-5}\}, \{v_{i}\})$ ($(P^L_{n+5} + P^R_n, \{v_0, v'_{i-5}\}, \{v'_0, v_{i}\})$  resp.), which is equivalent to $(P^L_{i} + P^R_{i-5}) + (P^L_{n-1-i} + P^R_{n-1-i})$ ($(P^L_{i} + P^R_{i-5}, \{v_0\}, \{v'_0\}) + (P^L_{n-1-i} + P^R_{n-1-i})$. Here, we have $(P^L_{n-1-i} + P^R_{n-1-i}) \equiv 0$ and $(P^L_{i} + P^R_{i-5})$ ($(P^L_{i} + P^R_{i-5}, \{v_0\}, \{v'_0\})$ resp.) satisfies the induction hypothesis as $n+4 \ge i \ge n \ge 6$ ($i \ge n \ge 7$ resp.) and therefore the score is at least one.
\end{itemize}

This strategy ensures that $Rs(P^L_{n+5} + P^R_n) \ge 1$ ($Rs(P^L_{n+5} + P^R_n, \{v_0\}, \{v'_0\}) \ge 1$ resp.).
\end{proof}

\subsubsection{Strategy for Right when Left starts}

When Left starts, the induction made in the previous proof cannot be applied. Indeed, from the position $(P^L_{n+5} + P^R_n, \{v_0\}, \{v'_0\} )$, Left can in one move make the position be $(P^L_{n+5} + P^R_n, \{v_0, v_{n+3}\}, \{v'_0\} )$ and no move of right can transform it into a position handled by the induction hypothesis. Therefore, another strategy is required. We will consider a strategy for Right that consists, for the leftmost vertices of both paths, in mimicking any move of Left on the other path, and that ensures some minimal properties on the moves played on the rightmost vertices. We introduce the following lemma to handle the rightmost vertices.

\begin{lemma} \label{strategy right}
Consider the graph $G = P^L_6 + \{v'_0\}$. Let $v_0$ be an extremity of $P^L_6$. In \mbinc{}, Right has a strategy, going second, such that Left claims either $v_0$ and $v'_0$ without any point, or at most one of $\{v_0, v'_0\}$ and she scores at most one point on $G$.
\end{lemma}

\begin{proof}

Let $G = P^L_6 + \{v_0'\}$. Recall that $v_0, \dots, v_5$ are the vertices of $P^L_6$. We will describe a strategy for Right playing second such that Left scores no point or such that she does not claim both $v_0$ and $v_0'$ with at most one point. 

\begin{itemize}
    \item If Left plays $v_0$, Right answers $v_1$,
    \begin{itemize}
        \item if Left plays $v_0'$, Right plays $v_3$ and pairs $v_4$ and $v_5$. Left cannot score a point. 
        \item If Left plays $v_2$ ($v_5$ resp.), Right plays $v_3$ ($v_4$ resp.) and pairs $(v_4, v_5)$ ($(v_2, v_3)$ resp.). This way, Left cannot score a point.
        \item If Left plays in $v_3$ ($v_4$ resp.), Right plays $v_4$ ($v_3$ resp.) and pairs $v_2$ ($v_5$ resp.) with $v_0'$. Either Left scores a point or claims both $v_0$ and $v_0'$. 
    \end{itemize}
    \item If Left plays $v_1$, Right answers $v_0$. He has claimed one of $(v_0, v_0')$. He then pairs $(v_2, v_3)$ and $(v_4, v_5)$. With this pairing, Left can score at most one point.
    \item If Left plays $v_2$, Right answers $v_3$, He then pairs ($v_0, v_1$ and $v_4, v_5$). The only one edge outside the pairing (and therefore that can be claimed by Left) is $v_1, v_2$ but with this pairing, Right then plays $v_0$ and claim one of $v_0, v_0'$. Otherwise, Left scores no point.
    \item If Left plays $v_3$ ($v_5$ resp.), Right answers $v_4$, he then pairs $(v_0, v_0')$ and $(v_1, v_2)$. This way, Left scores at most one point on the edge $(v_2, v_3)$ or $(v_0, v_1)$ but she cannot take both. And Right will be able to take one of $v_0$ or $v_0'$.
    \item If Left plays $v_4$, Right answers $v_3$
    \begin{itemize}
        \item If Left plays $v_0$, Right plays $v_1$ and pairs $v_0'$ with $v_5$. Either Left claims $v_0'$, and then by claiming $v_5$, Right ensures that Left scores no point, or Left claims $v_5$ and scores one point, but Right claims $v_0' \in \{v_0, v_0'\}$.
        \item If Left plays $v_1$ ($v_2$ resp.), Right plays $v_0$ and pairs $v_2$ ($v_1$ resp.) and $v_5$. By claiming one of them, Left scores one point but Right claims the second one, and therefore, Right ensures that Left scores only one point and does not claim both $v_0$ and $v_0'$.
        \item If Left plays $v_5$ ($v_0'$ resp.), Right plays $v_0$ and pairs $(v_1, v_2)$. Then, Left cannot score a second point (can score at most one point by playing $v_5$ resp.), and Right has already claimed one of $v_0, v_0'$.
    \end{itemize}
    \item If Left plays $v_0'$, Right answers $v_0$. He has already claimed one of $v_0, v_0'$, and the remaining graph is equivalent to  $P^L_5$ for which we already know that Left gets at most $1$ when she starts.
\end{itemize}
\end{proof}

\begin{lemma}\label{score Ls path}
Let $n \ge 1$ be an integer. In \mbinc{}, we have $ Ls(P^L_{n+5} + P^R_n) \le 1$.

Let $n \ge 2$ be an integer. In \mbinc{}, we have $ Ls(P^L_{n+5} + P^R_n, \{v_0\}, \{v'_0\}) \le 1$.
\end{lemma}

\begin{proof}
The proof below holds for the two cases, i.e. if the vertices $v_0$ and $v'_0$ are already colored or not.

Recall that $v_0, \dots, v_{n+4}$ are the vertices of $P^L_{n+5}$ and $v'_0, \dots, v'_{n-1}$ are the vertices of $P^R_n$.
We provide here a strategy for Right to ensure that the score is at most $1$ as follows:

\begin{itemize}
    \item If Left plays a vertex in a pair ($v_i$, $v'_i$) with $0 \le i \le n-2$, Right answers the second vertex of this pair.
    \item If Left plays another vertex, Right follows the strategy of Lemma~\ref{strategy right} with $P_6^L = \{v_0=v_{n-1},\ldots,v_{n+4}\}$ and $v_0' = v'_{n-1}$.
\end{itemize}

According to this strategy, Right ensures that Left scores the same number of points as him on the subgraph induced by the vertices $v_i, v'_i$ with $0 \le i \le n-2$. On the rest of the graph, from Lemma~\ref{strategy right}, either Left takes the two vertices $v_0,v'_0$ and gets no point, which can yield her overall at most one point with the edge $(v_{n-2},v_{n-1})$ of $P^L_{n+5}$. Otherwise, she takes $v'_{0}$ or the extremity $v_{0}$  of the $P_6^L$ and scores one point. In this case, if this extremity corresponds to $v'_{n-1}$ of $P^R_n$ she does not score a second point, and if this extremity is $v_{n-1}$, she can score a point if she also takes $v_{n-2}$. But in this case, Right has claimed both $v'_{n-2}$ by the pairing strategy and $v'_{n-1}$ as he has also claimed the other extremity. Thus, Right also scores one point. Finally, Right ensures that the score is at most $1$ with this strategy, and we have $Ls(P^L_{n+5} + P^R_n) \le 1$.
\end{proof}

\subsubsection{Proof of Theorem~\ref{equivalence paths} and score on paths}

Now we can prove Theorem~\ref{equivalence paths}.

\begin{proof}
By symmetry, as $P^L_n = -P^R_n$ for any $n$, we only need to prove the result for $P^L_n$.

As our game is in Milnor's universe, according to Lemma~\ref{lemma equiv}, it is sufficient to prove that $P^L_{n+5} - P^L_n - 1 \equiv 0$ ($(P^L_{n+5} - P^L_n, \{v_0\}, \{v_0'\}) - 1 \equiv 0$ resp.), i.e. $Ls(P^L_{n+5} + P^R_n) = Rs(P^L_{n+5} + P^R_n) = 1$ ($Ls(P^L_{n+5} + P^R_n, \{v_0\}, \{v_0'\}) = Rs(P^L_{n+5} + P^R_n, \{v_0\}, \{v_0'\}) = 1$ resp.). 

As the game is nonzugzwang, and according to Lemma~\ref{strategy left} and Lemma~\ref{score Ls path}, we have proven $1 \ge Ls(P^L_{n+5} + P^R_n) \ge Rs(P^L_{n+5} + P^R_n) \ge 1$ ($1 \ge Ls(P^L_{n+5} + P^R_n, \{v_0\}, \{v_0'\}) \ge Rs(P^L_{n+5} + P^R_n, \{v_0\}, \{v_0'\}) \ge 1$ resp.), which corrresponds to the desired result.
\end{proof}

From Theorem~\ref{equivalence paths}, and since the score on small paths is provided by Figure~\ref{fig:score_path}, the score of any path can be computed as follows:

\begin{corollary}
Let $n \ge 1$ be an integer. Denote by $n = 5q + r$ with $q$ and $r$ the quotient and the rest of $n$ divided by $5$. In \mbinc{}, on the one hand, we have $Ls(P^L_n) = -Rs(P^R_n) = q$ if $0 \le r \le 2$, and $Ls(P^L_n) = -Rs(P^R_n) = q +1$ if $3 \le r \le 4$. On the other hand, we have $Rs(P^L_n) = -Ls(P^R_n) = q-1$ if $r = 0$, $Rs(P^L_n) = -Ls(P^R_n) = q$ if $1 \le r \le 4$.
\end{corollary}

\subsection{Union of paths and cycles}

We will denote cycles as follows:
\begin{itemize}
\item $C^L_n$: cycle of length $n$ where all the edges are colored blue.

\item $C^R_n$: cycle of length $n$ where all the edges are colored red.
\end{itemize}

Now that the equivalences of paths are known, union of paths can easily be reduced to union of paths of order at most $5$. Yet, to deal with such unions, it is not sufficient in general to compute the score on them. The problem can be solved by considering new equivalences between small paths.

\begin{lemma}\label{small equivalences}
In \mbinc{}, we have the following equivalences:
\begin{align}
    P^L_1 \equiv P^L_2 \equiv 0 \\
    2P^L_3 \equiv 1 \\
    P^L_4 \equiv P^L_3 \\
    2 P^L_5 + P^L_3 \equiv 2
\end{align}
\end{lemma}

\begin{proof}
Recall that given a graph $G$ and an integer $k$, in order to prove that $G \equiv k$, it is sufficient to prove $k \ge Ls(G)$ and $Rs(G) \ge k$.

\begin{enumerate}
    \item We have $Ls(P^L_1) = Rs(P^L_1) = 0$ and $Ls(P^L_2) = Rs(P^L_2) = 0$ as in both games no edges are taken by a player. This proves, by Lemma~\ref{lemma equiv}, that $P^L_1 = P^L_2 = 0$
    \item We prove $Ls(2P^L_3) = Rs(2P^L_3) = 1$. To do that, we just need to prove $Ls(2P^L_3) \le 1$ and $Rs(2P^L_3)\ge 1$. 
    \begin{itemize}
        \item Suppose Left starts. If she plays in one path $P^L_3$, Right claims the middle vertex of the other path and then plays at least one vertex in the $P^L_3$ where Left started. This way, Left scores at most one.
        \item Suppose Right starts. He plays in one path $P^L_3$. By going first in the second path, Left can score one by playing the middle vertex and after that at least one of its two neighbors. 
    \end{itemize}
    \item As $- P^L_3 = P^R_3$, we will prove $P^L_4 + P^R_3 = 0$. Denote by $(v_0, v_1, v_2, v_3)$ the vertices of $P^L_4$ and by $(v'_0, v'_1, v'_2)$ the vertices of $P^R_3$
    \begin{itemize}
        \item Suppose Left starts. If she plays in $P^R_3$, Right plays $v_1$ and pairs $(v_2, v_3)$ to ensure that Left cannot score an edge. If Left plays $v_0$ or $v_1$ ($v_2$ or $v_3$ resp.), Right plays $v_2$($v_1$ resp.) and pairs $(v'_0, v'_2)$ and $v'_1$ with the available vertex in $\{v_0, v_1\}$ (in $\{v_2, v_3\}$ resp.). This way, Left and Right scores the same number of edges and this proves $Ls(P^R_3 + P^L_4) \le 0$
        \item Suppose Right starts. Left considers the pairing $(v_0, v'_0)$, $(v_1, v'_1)$, $(v_2,v'_2)$. This way, any point scored by Right is scored by Left. Therefore $Rs(P^L_4 + P^R_3) \ge 0$.
    \end{itemize}
    \item Let $G = 2P^L_5 + P^L_3$. Denote by $v_0, \dots, v_4$ and $v_0', \dots, v_4'$ the vertices of the two copies of $P^L_5$ and by $(u_0, u_1, u_2)$ the vertices of $P^L_3$. Let first prove $Rs(G) \ge 2$. Up to consider only $3$ vertices of one copy of $P^L_5$, we can suppose that the first move of Right is in a $P^L_5$ and we will prove that Left scores $2$ on $P^L_5 + P^L_3$. Suppose Right has played a vertex $v_i'$ with $0 \le i \le 4$. Left plays $v_2$ and continues as follows:
    \begin{itemize}
        \item If Right plays $v_0$ or $v_1$ ($v_3$ or $v_4$ resp.), Left plays $v_3$ ($v_1$ resp.) and pairs $(v_4, u_1)$ ($(v_0, u_1)$ resp.) and $(u_0, u_2)$.
        \item If Right plays $u_0, u_1$ or $u_2$, Left plays $v_1$ and pairs $(v_0, v_3)$.
    \end{itemize}
    In both cases, Left scores at least two points.
    Now we prove $Ls(G) \le 2$. After the first move of Left, at least one of the two copies of $P^L_5$ has its $5$ vertices available. Suppose it is $v_0', \dots, v_4'$. Right plays $v_2'$ and pairs $(v_0', v_1')$ and $(v_3', v_4')$, ensuring Left won't score any point on this copy of $P^L_5$. Left plays a second move:
    \begin{itemize}
        \item If $v_2$ has not been played yet, Right plays $v_2$. Left plays a third move. If the three moves of Left are in $P^L_3$, Right pairs $(v_0, v_1)$ and $(v_3, v_4)$, ensuring Left does not score any other point. If at least one of them is not in $P_3$, Right plays any vertex of $P_3$, and know that at least one vertex of $(v_0, v_1, v_3, v_4, u_0, u_1, u_2)$ will be available for his next move. Thus, Left cannot score more than two points on the rest of them.
        \item If Left has played $v_2$, at least one of $v_1$ or $v_3$ is available. Right plays it. By symmetry, suppose it is $v_1$. After the next move of Left, at least one of $v_3, v_4, u_1$ will be available. Right plays it, ensuring again that Left cannot score more than $2$.
    \end{itemize}
\end{enumerate}
\end{proof}

We can now state the equivalence theorem for union of paths.

\begin{corollary}\label{corollary equiv path}
Let $P_1, \dots, P_N$ be paths of lengths $n_1, \dots, n_N$.

Let $q_1, \dots, q_N$ be positive integers and $
1\le r_1, \dots, r_N \le 5$ be integers such that for any $1 \le i \le N$, we have $n_i = 5q_i + r_i$. 

Denote for $1 \le i \le 5$ by  $N_i$ the number of $r_j$ equal to $i$. In \mbinc{}, we have:
$$ \underset{i = 1}{\overset{N}{\sum}} P^L_{5 q_i + r_i} \equiv  \underset{i = 1}{\overset{N}{\sum}} q_i + \left \lfloor \frac{N_3 + N_4}{2} \right \rfloor + 3\left \lfloor \frac{N_5}{4} \right \rfloor + (N_3 + N_4 \mod 2) P_3 + (N_5 \mod 4) P_5 $$

Therefore, $Ls(\underset{i=1}{\overset{N} \sum} P_i)$ and $Rs(\underset{i=1}{\overset{N} \sum} P_i)$ are computable in linear time.
\end{corollary}

\begin{proof}
By Theorem~\ref{equivalence paths}, any path $P^L_{5 q_i + r_i}$ is equivalent to $q_i + P_{r_i}$. Then, by Lemma~\ref{small equivalences}, we have $P_3 \equiv P_4$, $2P_3 \equiv 1$, and $2 (2P_5 + P_3) \equiv 4 P_5 + 2 P_3 \equiv 4 P_5 + 1 \equiv 4$. Thus $4P_5 \equiv 3$. Note that these computations are possible thanks to Milnor's universe.
\end{proof}

Note that we consider $1 \le r_i \le 5$ and not $0 \le r_i \le 4$, so $q_i$ and $r_i$ are not exactly the quotient and the rest of the size of the path by $5$.

\begin{corollary}
Let $n \ge 1$. In \mbinc{}, there exists a linear time algorithm to compute $Ls(C^L_n)$ and $Rs(C^L_n)$.
\end{corollary}

\begin{proof}
First, note that $Rs(C^L_n) = Ls( P^L_{n-1})$.

To compute $Ls(C^L_n)$, note that all the vertices are symmetric. Therefore, we can suppose that Left first plays any of them. The next move of Right will make the graph equivalent to $(P^L_{k}, \{v_0\}, \emptyset) + (P^L_{k'}, \{v_0'\}, \emptyset)$ with $v_0, v_0'$ extremities of $P^L_k$ and $P^L_{k'}$ and with $k + k' = n$. The score on these graphs can be computed in linear time by using Corollary~\ref{corollary equiv path}, and therefore, $Ls(C^L_n)$ too as, by Theorem~\ref{equivalence paths}, at most $5$ values are to be considered for the pair $(k,k')$ according to the equivalences.
\end{proof}

\section{Perspectives}

In this paper, we introduced positional scoring games in a general framework and then focused on \inc, which corresponds to the case of 2-uniform hypergraphs. 
To conclude this paper, we list some relevant open problems.

\begin{itemize}
    \item We have solved \mbinc{} on union of paths using game equivalences. Next step would be to study trees.
    \item What is the complexity of \mbinc{} when restricted to the class of cographs? Equivalent vertices have an important role and can be easily simplified. This could be a starting point for the study of cographs.
    \item We proved that \mbinc{} is fixed-parameter tractable using the neighborhood diversity. It would be interesting to find other parameters for which the problem is FPT. For example, is it FPT parameterized by the score?
    \item The same question applies when considering general hypergraphs. The answer is negative for 6-uniform hypergraphs as it is {\sf PSPACE} even for $k=1$. What about $3$-uniform hypergraphs?
    Since \MBPG{} is polynomial for $3$-uniform hypergraphs \cite{galliot2022}, the question makes sense.
    \item We have proved that \MMSPG{} is {\sf PSPACE}-complete even for 3-uniform 
    hypergraphs but provided a linear algorithm for $2$-uniform hypergraphs. It might be interesting to look at particular $3$-uniform hypergraphs. For example, is it possible to compute the score in the scoring version of the {\sc Triangle Game} (where players choose edges of a graph and try to construct triangles)? The hypergraph of this game has the particularity to be linear (hyperedges cross on at most one vertex). A more general question would be to find the complexity of $\mminc{}$ on linear $3$-uniform hypergraphs.
    \item In Section \ref{sec:milnor}, we have introduced {\em partisan scoring positional games} to include the two conventions of scoring positional games in a more general definition. Maker-Maker convention corresponds to games with only green hyperedges whereas Maker-Breaker convention corresponds to games with only blue edges. It would be interesting to consider games with both red and blue edges but no green edge.
\end{itemize}

\bibliography{Incidence} 
\bibliographystyle{alpha}

\end{document}